\documentclass[journal,onecolumn]{IEEEtran}
\usepackage[numbers]{natbib}
\usepackage[T1]{fontenc}
\usepackage[utf8]{inputenc}
\usepackage[english]{babel}
\usepackage{amsmath, amsthm, amssymb}
\usepackage{hyperref}
\usepackage{url}
\usepackage{physics} 
\usepackage{graphicx}
\usepackage{mathtools}
\usepackage{thmtools}
\usepackage{etoolbox}
\usepackage{amssymb,setspace,amsthm,calrsfs}
\usepackage{setspace}
\usepackage{pgfplotstable,pstool,psfrag}
\usepackage{float}
\usepackage{tabularx}
\usepackage{ltablex}
\usepackage{booktabs}
\usepackage{algorithm2e}
\usepackage{picins}
\usepackage{tabularx}
\usepackage{supertabular}
\usepackage{stmaryrd}
\usepackage{scalerel}
\usepackage{tikz-cd}
\usepackage{color}
\usepackage{ulem}
\usepackage[bottom]{footmisc}
\newtheorem{thm}{Theorem}[section]

\newtheorem{dfn}[thm]{Definition}
\theoremstyle{remark}
\newtheorem{rem}[thm]{Remark}

\newcommand{\R}{\mathbb{R}}
\newcommand{\N}{\mathbb{N}}

\newcommand{\probability}{\mathbb{P}}

\def\defeq{\mathrel{\mathop:}=}

\pgfplotsset{compat=1.6}

\makeatletter
\let\mytagform@=\tagform@
\def\tagform@#1{\maketag@@@{\color{green!55!black}(#1)}}
\makeatother

\newif\iffinalfonts
\iffinalfonts
\else
\newcommand{\bbPi}{\Pi}

\fi

\title{Predictability and Fairness in the Interconnections of Ensembles with Applications to Two-Sided Markets}
\author{Wynita M. Griggs, Ramen Ghosh, Jakub Mare\v{c}ek, and Robert N. Shorten%
\thanks{This work was supported in part by the Science Foundation Ireland under Grant 16/IA/4610. Jakub acknowledges support of the OP RDE funded project CZ.02.1.01/0.0/0.0/16 019/0000765 ``Research Center for Informatics''.}
\thanks{W. M. Griggs is with the Department of Civil Engineering and the Department of Electrical and Computer Systems Engineering, Monash University, Clayton, Victoria, 3800, Australia.}
\thanks{R. Ghosh is with the School of Electrical and Electronic Engineering, University College Dublin, Ireland.}
\thanks{J. Marecek is with the Czech Technical Univeristy in Prague, the Czech Republic.}
\thanks{R. N. Shorten is with the the Dyson School of Design Engineering, Imperial College London, South Kingston, UK.}
}

\newtheorem{prop}[thm]{Proposition}

\begin{document}

\maketitle

\begin{abstract}
There has been much recent interest in two-sided markets and dynamics thereof. In a rather a general discrete-time feedback model, which we show conditions that assure that for each agent, there exists the limit of a long-run average allocation of a resource to the agent, which is independent of any initial conditions. We call this property the unique ergodicity.
  Our model encompasses two-sided markets and more complicated interconnections of workers and customers, such as in a supply chain. It allows for non-linearity of the response functions of market participants. Finally, it allows for uncertainty in the response of market participants by considering a set of the possible responses to either price or other signals and a measure to sample from these. 
\end{abstract}

\section{Introduction}
Motivated by the success of the business models of Uber Technologies, Inc.\ or Upwork Global Inc.,
and following Jean Tirole's pioneering research \citep{Rochet2003,Rochet2004(1),Rochet2004(2),Rochet2006},
there has been much recent work on two-sided and multi-sided markets (e.g., by \cite{Weyl2010,Hagiu2015,Evans2016,Parker2016(1)} and most recently by \cite{ashlagi2020clearing,Arnosti2021,Benjaafar2021,Liu2021,kanoria2018convergence,mohlmann2020algorithmic}).
One would like to analyse the long-run properties of such markets, including their stability and fairness.

We consider a discrete-time feedback model of the interconnection of ensembles, which:
\begin{itemize}
\item encompasses two-sided markets and more complicated interconnections in multi-sided markets,
\item allows for the nonlinearity of the response functions of market  participants,
\item and allows for uncertainty in market participants' responses by considering a set of possible responses to either price or other signals, and a measure to sample from these.
\end{itemize}
Our model applies even to settings where the organizing entity divides resources among agents, based on the information reported by the agents without payments, which is essential within ``artificial intelligence for social good''. 
While one can leverage auditing mechanisms to maximize utility in repeated allocation problems where payments are not possible \citep{Lundy2019}, our results can be seen as conditions on the information exchange, as well as prices, that allow for specific desirable properties.  
Our results concern specific desirable properties of such models and certain testable conditions assuring these properties \citep{Fioravanti2019}.
Informally speaking, we call a feedback system uniquely ergodic, when for every agent $i$, there exists a limit of a long-term average allocation of a resource to the agent, independent of any initial conditions.
In turn, the notion of unique ergodicity underlies a natural notion of fairness, distinct from other popular notions of fairness \citep{Bateni2016}, where the limit coincides across all agents or market participants.
A necessary condition for the feedback system to be uniquely ergodic is to be contractive on average. In a sense, we formalise below. 
\subsection{An Application: Two-sided Markets}
Critical applications lie within two-sided markets \citep{Rochet2004(1), Rochet2004(2), Rochet2006, Lobel2020}, which model online labour platforms that enable the interaction between the two sides: customers submitting jobs over time and workers performing jobs. 
For example, the ride-hailing systems of Uber, Lyft, or Didi Chuxing can be modelled as a two-sided market with direct connections. In this case, workers would be Uber drivers, and jobs can be the rides of customers. Note that jobs are assumed to be independent of each other, although one customer might offer multiple jobs.

See Figure~\ref{fig:system_ii} for an overview of our model of such a two-sided market:
controller ${\mathcal C^1}$ suggests prices $\pi^1$ (based on the distance travelled and the so-called driver surge pricing \citep{Chen2016,Castillo2017,Cachon2017,Castillo2020,Garg2020} in Uber) to customers ${\mathcal S}^1_i$, whose requests $y^{1}_i(k)$ for jobs (rides) at time $k$ are based on some internal state of each customer $i$ at time $k$, $x^{1}_i(k)$, which are not directly observable. A controller ${\mathcal C^2}$ for the other side of the market matches the jobs (rides) to workers (driver-partners) ${\mathcal S}^2_i$ whose state $x^{2}_i(k)$ at time $k$ may be partially observable (e.g., availability, position) and partially not observable (e.g., appetite for further work that day). 
Usually (e.g., \cite{Simonetto2019,Araman2019,Aouad2020,Ozkan2020,Yan2020}), ${\mathcal C^2}$ is implemented using an on-line matching algorithm. Its matches are provided to the workers (drivers), whose total number $y^{1}(k)$ of accepted matches is then filtered to obtain the proportion of the empty cars on the road, which is then the input into the controller ${\mathcal C^1}$ that suggests prices (with driver surge pricing implemented, if there are too few empty cars). 
\subsection{Another Application: Multi-sided Platforms}
\label{app:multi-sided}
Further important applications lie within multi-sided platforms \citep{Weyl2010,Hagiu2015,Evans2016} and networked markets \citep{Parker2016(1)}.
To continue our ride-hailing example, one could see the ride-hailing system of Uber as a multi-sided market if one also considers ``fleet partners'' (who are intermediaries for car manufacturers and car leasing providers) and taxi operators. Indeed, Uber's Vehicle Solutions Program is a platform for fleet partners to offer their vehicles to driver-partners. Drivers, in turn, sometimes happen to work also as licensed taxi drivers.  

Similarly, Google's Android ecosystem and  Microsoft's Windows ecosystem are sometimes seen \citep{Hagiu2015, Parker2016(2)} as three-sided platforms connecting consumers, software providers, and hardware providers.
While some of the incentives offered to independent software providers (e.g., free hardware samples, no-cost licenses of development tools) are not being adjusted in real-time, others (e.g., promotions for their apps) are. The details of the control mechanisms have not been made public in this case. 
\subsection{Related Work}
Within Economics, Jean Tirole's pioneering research \citep{Rochet2003,Rochet2004(1),Rochet2004(2),Rochet2006} on two-sided markets has a substantial 
following \citep{Lobel2020}.
Especially multi-sided platforms \citep{Weyl2010,Hagiu2015,Evans2016} and networked markets \citep{Parker2016(2)}
became widespread both in economic theory and the real world. 

Within Economics and Computation, much attention has been focused on ride-hailing. Dynamic pricing cannot outperform the optimal static policy in terms of throughput and revenue when workers cannot reject jobs. 
It is known \citep{Garg2020} that additive increases provide an incentive-compatible pricing mechanism,
while certain complications arise \citep{Ma2020} from the spatio-temporal nature of the problem. 

Motivated by Uber's admission \citep{Cook2018} that its prices are unfair to female and ethnic minority drivers, there has recently been some interest in studying fairness in related settings \citep{Suhr2019, Cohen2019, Jung2020, Freeman2020}, often using 
concept of fairness \citep{Chouldechova2018,Mouzannar2019} 
articulated in machine learning. 

Related work is also described in \cite{Agents07} using ideas from multi-agent systems; and in \cite{Blondel05, Nedic09} from the perspective of 
distributed optimization, in both cases exploiting the relationship between consensus, utility maximization, and fairness.
The concept of fairness is also strongly connected to {\em ergodicity} as ergodic dynamic behaviour implies several salient features 
that is necessary for fairness \citep{Mezic11}.
Finally, we note that our results in this paper are established building upon the seminal work on {\em iterated random functions}
\citep{Elton1987,Barnsley1988(1),Barnsley1989}. Although not well known, iterated random functions (IRF) are a class of discrete-time Markov processes with sufficiently rich background results. These types of systems provide a natural framework for modelling classes of multi-agent systems, and in this setting, a wealth of known and established results can be applied to analyze such systems \citep{Elton1987, Barnsley1988(2), Barnsley1989, Barnsley2013, Stenflo2001(1), Szarek2003(1), Steinsaltz1999, Walkden2007, Barany2015, Diaconis1999, Iosifescu2009, Stenflo2012(s)}. 
\subsection{Contributions and Paper Organization}
\begin{figure*}[htbp]
\centering
\includegraphics[width=0.75\columnwidth]{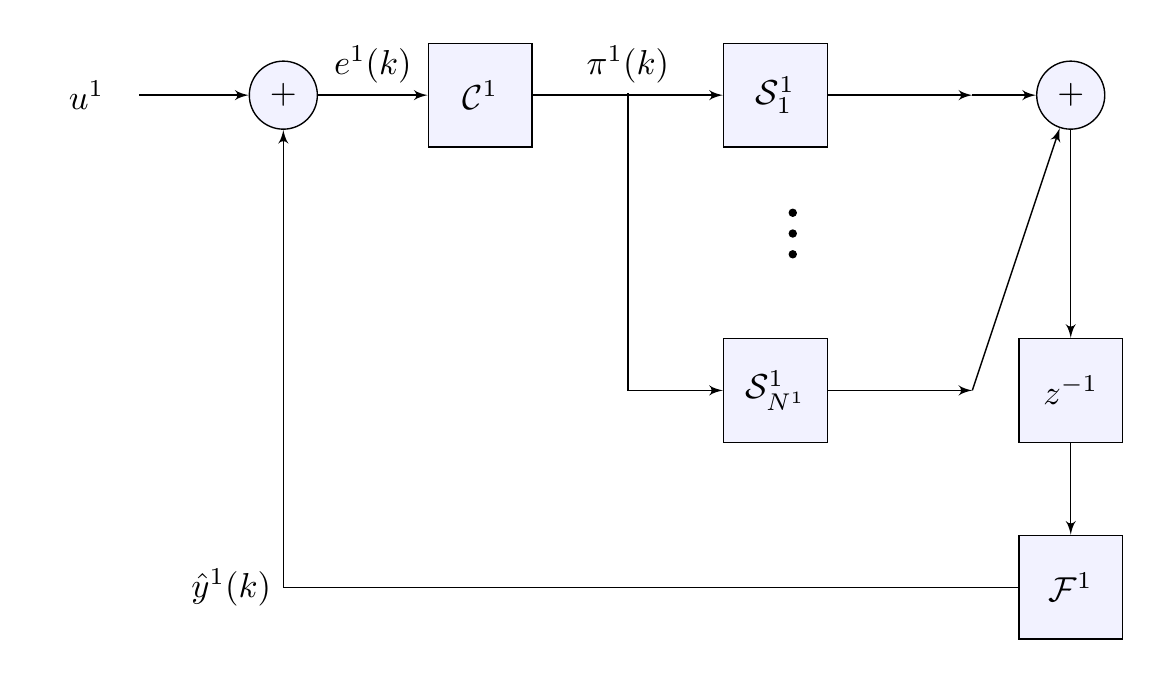}
\caption{The feedback model of \cite{Fioravanti2019}, depicted utilizing our notation.}\label{fig:system}
\end{figure*}

\begin{figure*}[htbp]
\centering
\includegraphics[width=0.8\columnwidth]{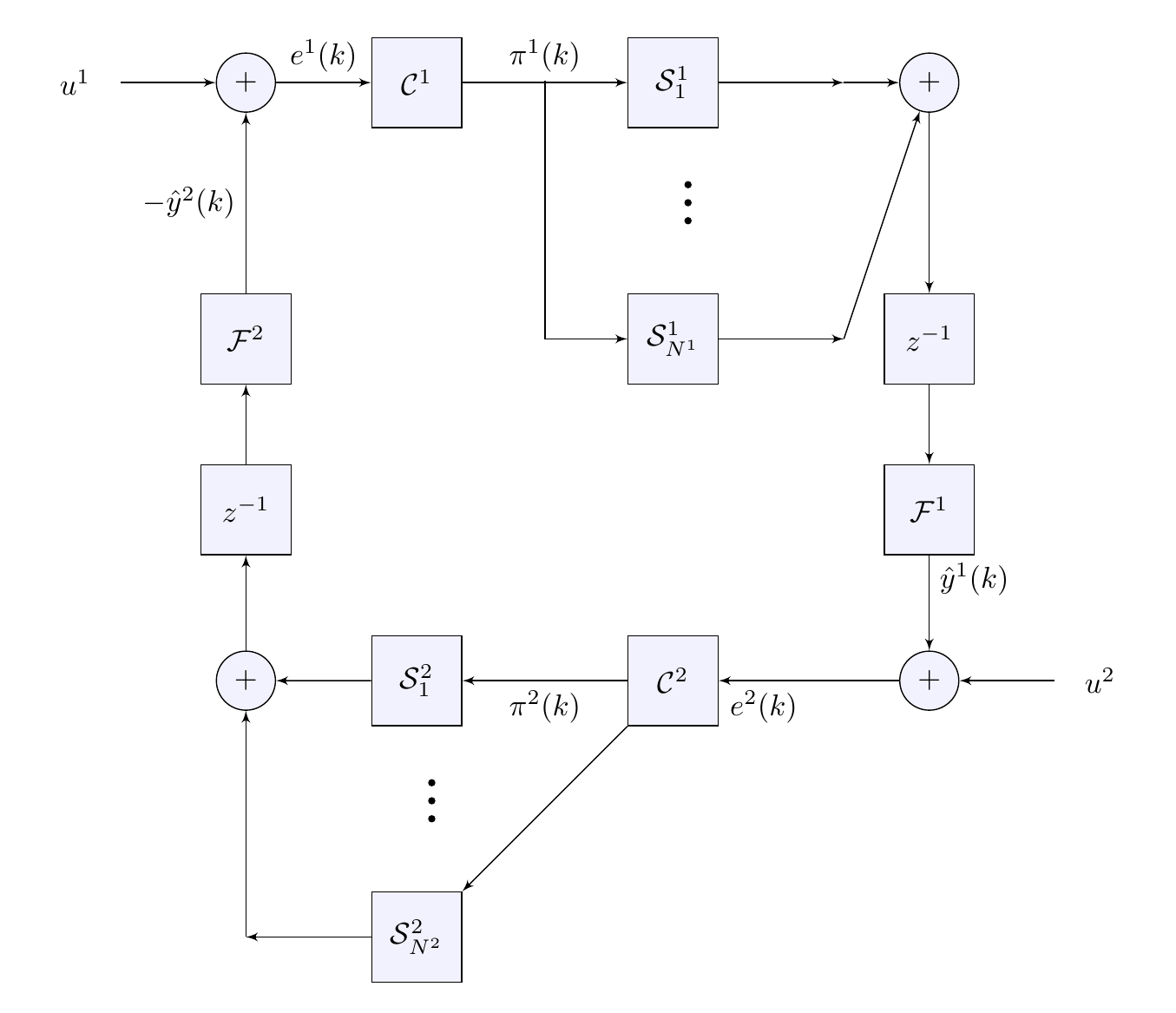}
\caption{Our feedback model of a two-sided market.}\label{fig:system_ii}
\end{figure*}

For these -- and many other -- applications, we present: 
\begin{itemize}
    \item a novel model of two-sided markets and extensions towards more general interconnections of ensembles of agents,
    \item a notion of unique ergodicity, which can be seen as individual-level predictability of the outcomes in a certain closed-loop sense,
    \item and conditions ensuring unique ergodicity.
\end{itemize}
Our model is based on our earlier work on ensemble control  \citep{Fioravanti2019}, where all agents respond to the prices a central organizing entity sets or signals it provides. We extend this to multiple populations of agents in a multi-sided market and more complicated interconnections. 
The differences are illustrated in Figures \ref{fig:system} and \ref{fig:system_ii}: Figure \ref{fig:system} is the model used in our earlier work; Figure \ref{fig:system_ii}, the present model.


Our results are based on {\em iterated random functions}
\citep{Elton1987,Barnsley1988(1),Barnsley1989}. \emph{Iterated random functions} (IRF) are not prominent in the Economics community. Nevertheless, we will see that by using IRF, two-sided markets can be modelled and analysed naturally. That makes it possible to establish strong stability guarantees.

The paper is organised as follows. In Section~\ref{sec:preliminaries}, we describe our model and provide the necessary mathematical background on {\em discrete time Markov processes}, {\em iterated random functions}, and briefly note the idea of {\em coupling of invariant measures} to point out a necessary condition for {\em ergodicity}. Our main results are established in Section \ref{sec:inter-two-ens} on the statistical stability of a two-sided market, which is modelled by the closed-loop system as in Figure \ref{fig:system_ii}, where in Theorem \ref{thm:ITE-NLS-LCF} we state the result assuming that the agents' states evolve as a {\em linear time-invariant dynamical system}; in Theorem \ref{thm:intrcon-twoensmbl-nonlinrstate-lincontr&filtr} the states are realized by a {\em nonlinear iterated random function}; and in Theorem \ref{thm:discrete-range-space} we assume that the agent's state has discrete-range space. In Section \ref{sec:inter-lrgscl-ensmbl} we extend our results obtained in Section \ref{sec:inter-two-ens} for the interconnection of $N$ number of ensembles of systems.
And finally in Section \ref{sec:conclsn}, we highlight some directions for future research.

\section{Preliminaries, Notation, and the Description of our Model}\label{sec:preliminaries}
In this article, we intend to extend and generalize our results from \cite{Fioravanti2019} about the closed-loop,
discrete-time dynamical systems of the form depicted in Figure \ref{fig:system}, which consist of a number $N \in \N$ denoting the total number of agents, $\mathcal F$ denoting the filter, and a central controller  $\mathcal{C}$ that produces a signal $\pi(k)$ at time $k$. In response, the agents, modelled by the systems $\mathcal{S}_1,\mathcal{S}_2$, \ldots, $\mathcal{S}_{N}$, modify their use of the resource.

Let $\left(\mathbb{X}, {\mathcal B}\right)$ denote a {\em measurable state-space}. In most of our discussions, $\mathbb X$ is the Eucleadian space $\mathbb R^n$ or a closed subset of it, let $\rho$ be the usual Eucleadian metric, and $\mathcal B$ is a Borel $\sigma$-algebra on it. $\mathcal M\left(\mathbb{X}\right)$ denotes the space of all {\em Borel probability measures} on $\mathbb{X}$.
For $i=1,2,\dots, N$, let $x_i(k)$ denote the internal state of the agents in the network. In particular, for all $i=1,2,\dots,N$, $x_i(k)$ is a random variable. Furthermore, the resource of agent $i$ at time $k$ is denoted by $y_i(k)$ and is modeled as the output of system ${\mathcal S}_i$. 
Here $y_i(k)$ and $\pi(k)$ are scalars, but this is easily generalized. The randomness in the system can be a consequence of the inherent randomness in the response of user $i$ to the control signal $\pi(k)$, or the response to a control signal that is intentionally randomized \citep{Arieh13, Arieh14,marevcek2015signaling}.
It is easy to see that the total resource utilisation
\begin{align}\label{eq:agg-res-uti}
y(k) \defeq \sum\limits_{i = 1}^N y_i(k)
\end{align}   
is also a random variable. The controller will not have information of any of the $x_i(k)$, $y_i(k)$ or $y(k)$, but it will always have access to $e(k)$ which is an error signal and expressed as follows 
\begin{align*}
e(k):=\hat y(k)- r,    
\end{align*}
where
$\hat y(k)$ is the output of a filter $\mathcal{F}$, and $r$ is a {\em desired value} of $y(k)$.

Let $\bbPi$ denote the set of all {\em possible broadcast control signals}. Let the {\em private state} of the controller be denoted by $x_c(k) \in \R^{n_c}$, and at time instant $k$, it is intended to modulate the system in Figure \ref{fig:system} by sending a signal $\pi(k) \in \bbPi \subseteq \R$. In the simplest situation, $\pi(k)$ is a function of $e(k)$ and $x_c(k)$, whose range is $\bbPi$.

For ensembles of {\em discrete} agents, the non-de\-ter\-min\-is\-tic agent-specific response to the feedback signal $\pi(k)\in \bbPi$ can be modelled
by agent-specific and signal-specific probability distributions over certain agent-specific sets of actions $\mathbb{A}_i =\{a_1,\ldots,a_L \}\subset \R^{n_i}$.
Furthermore, we denote $\mathbb{D}_i$ as the set of possible resource demands of agent $i$, where in the finite case 
$\mathbb{D}_i := \{ d_{i,1}, d_{i,2}, \ldots, d_{i,m_i}\}.$

For general purpose, we can assume that there are $w_i \in \N$  maps ${\mathcal W}_{ij}: \R^{n_i} \to \R^{n_i}$, $j=1,\ldots, w_i$
for agent $i$ and $h_i \in \N$ output maps ${\mathcal H}_{i\ell}: \R^{n_i} \to \mathbb{D}_i$, $\ell= 1,\ldots,h_i$ for all agent $i$.
The dynamic system evolves according to:
\begin{align}\label{eq:general}
x_i(k+1) & \in  \{{\mathcal W}_{ij}(x_i(k)) \;\vert\; j = 1, \ldots, w_i\},\\
y_i(k) & \in  \{{\mathcal H}_{i\ell}(x_i(k)) \;\vert\; \ell = 1, \ldots, h_i\},
\end{align}
where each agent $i$'s response at time $k$ is determined according to the
functions $p_{ij} : \bbPi \to [0,1]$, $j=1,\ldots,w_i$, respectively
$p'_{i\ell} : \bbPi \to [0,1]$, $\ell=1,\ldots,h_i$. Specifically, 
\begin{subequations} \label{eq:problaws}
\begin{align}
&\mathbb{P}\big( x_i(k+1) = {\mathcal W}_{ij}(x_i(k)) \big) = p_{ij}(\pi(k)),\\
&\mathbb{P}\big( y_i(k) = {\mathcal H}_{i\ell}(x_i(k)) \big) = p'_{i\ell}(\pi(k)).
\label{eq:signals}
\intertext{Additionally, it also holds that}
&\sum_{j=1}^{w_i} p_{ij}(\pi) = \sum_{\ell=1}^{h_i} p'_{i\ell}(\pi) = 1.
\end{align}
\end{subequations}
We assume that the random variables $\{ x_i(k+1) \mid  i = 1,\ldots,N \}$ are conditioned on $\{ x_i(k)  \}, \pi(k)$, but independent. 
The  overall system can be modelled by the operator $P$ as $P: \mathbb{X} \times \bbPi \to \mathcal M(\mathbb{X})$. In order to reason about the evolution of the state, consider the space $\mathbb{X}^{\infty}$, and we introduce the space of probability measures on $\mathcal M(\mathbb{X}^{\infty})$ with the product $\sigma$-algebra.
The set-up resembles very closely that of {\em iterated random functions} in \cite{Elton1987,Barnsley1988(1),Barnsley1989}.
Iterated random functions are a class of stochastic dynamical systems, for which strong stability and convergence
results exist on compact or complete and separable metric spaces as a system's state space; see \cite{Barnsley1988(1),Diaconis1999, Steinsaltz1999,Stenflo1998,Stenflo1999,Stenflo2001(1),Stenflo2001(2),Stenflo2012(s),Szarek1997,Szarek1999,Szarek2000(1),Szarek2000(2),Szarek2000(3), Szarek2001, Szarek2003(1),Szarek2003(2),Ramen2019}. We now introduce some further concepts and results on iterated function systems.
\begin{dfn}\label{dfn:ifswsdp}
Let $\{f_i\}_{i=1}^{m}$ be continuous self-transformations on $\mathbb X$, and $\{p_i(x)\}_{i=1}^{m}$ be probability functions on $\mathbb X$, such that 
\begin{align*}
p_i(x):\mathbb X \to [0,1] \text{ for all } i\in [1,m], \text{ and } \sum\limits_{i=1}^{m} p_i(x)=1 \text{ for all } x\in \mathbb{X}.
\end{align*}
The pair of sequences
$
\left( f_1(x), f_2(x), \dots, f_m(x); p_1(x), p_2(x),\dots, p_m(x)\right)
$
is called an iterated random function with state(place)-dependent probabilities.
\end{dfn}
Any discrete-time Markov chain can be generated by an \emph{iterated random function with probabilities} (see \cite[Section 1.1]{Kifer2012} or \cite[Page 228]{Bhattacharya2009}), although such a representation is not unique (see \cite{Stenflo1999}).
Informally, the corresponding discrete-time Markov process on $\mathbb X$ evolves as follows:
choose an initial point $x_0\in \mathcal X$. Select an integer from the set $[1,m]$ in such a way that the probability
of choosing $\sigma$ is $p_{\sigma}(x_0)$, $\sigma\in [1,m]$. When the number $\sigma_0$ is drawn, define
$
x_1= f_{\sigma_{0}}(x_0).    
$
Having $x_1$, we select $\sigma_1$ according to the distribution 
$
p_1(x_1), p_2(x_1), \dots, p_m(x_1),    
$
and we define 
$
x_2= f_{\sigma_1}(x_1),
$
and so on.

Let us denote $P_n$ for $n=0,1,2,\dots$, the distribution of $x_n$; i.e.,
$
P_n(B)= \mathbb P(x_n \in  B) \text{ for some } B \in \mathcal B.
$
The above procedure can be formalised. For a given $x\in \mathcal X$ and a Borel subset $B\in \mathcal B $, we may easily show that the
transition operator for the given IFS is of the form
\begin{align}\label{eq:transii}
P(x, B):= \sum_{i=0}^{N} 1_{B}\left(w_{i}(x)\right) p_{i}(x).
\end{align}
$P(x,B)$ is the transition probability from $x$ to $B$, where $1_{B}$ denotes the characteristic function of $B$:
\[ 1_{B} := \begin{cases} 
1 & \text{ if } x\in B. \\
0 & \text{ if } x\in  B^c.
\end{cases}
\]
\begin{dfn}[Invariant Probability Measure]\label{dfn:invariant}
\citep{Barnsley1988(1),Szarek2000(3)} Suppose the initial condition $X(0)=x_0$ is distributed according to measure
$\mu$, (denoted by $\probability_\mu$). The random variable $X(k)$ is distributed
according to the measure $\mu_k$ and conditioned on an initial probability measure $\lambda$ with the following inductive relation:
\begin{align}\label{eq:measureiteration}
\mu_{k+1}(\mathbb{G})  := \int_{\mathbb X} P(x, \mathbb{G}) \, \mu_k(d x) \text{for all event } \mathbb{G} \in \mathcal{B}.
\end{align}
A measure $\mu^{\star}$ on $\mathbb X$ is said to be {\em invariant with respect to the Markov process} $\{ X(k) \}$ if $P\mu^{\star} = \mu^{\star}.$
\end{dfn}
\begin{dfn}
\emph{An invariant probability measure $\mu^{\star}$ is called
attractive}, if for every probability measure $\nu$, the sequence $\{
\mu_k \}$ defined by \eqref{eq:measureiteration} with initial
condition $\nu$ converges to $\mu^{\star}$ in distribution.
The existence of attractive invariant measures is
intricately linked to the ergodic properties of the system.
\end{dfn}
With this background, our general problem considered in this paper is modelled as a Markov chain with a state-space
representing all system components.
Our specific objective is to develop algorithms and systems to distribute the shared resource such that the following goals are achieved with a probability $1$. 
\begin{dfn}
Let $r>0$ be any upper bound for the utilization of the resource. Then we aim to have,
for all $k\in\N$,
\begin{align} \label{eq:feasi}
\sum_{i = 1}^N y_i(k) = y(k) \leq r.
\end{align}
In this article, we assume the upper bound to be a constant.
\end{dfn}
\begin{dfn}[Unique Ergodicity]\label{dfn:uniq-ergo}
We call a feedback system uniquely ergodic when, for every agent $i$, there exists a constant $\overline{r}_i$ such that
\begin{align}\label{eq:CLT}
\lim_{k\to \infty} \frac{1}{k+1} \sum_{j=0}^k y_i(j) = \overline{r}_i,
\end{align}
and $r_i$ does not depend on the initial condition of the agent's state.
\end{dfn}
Further to the above concept, one optional concern may be the idea of fairness, which could be conceptualized by saying that all the $\overline{r}_i$ coincide, thus the vector
$\overline{r} =
\left[\overline{r}_1 \ \ldots \ \overline{r}_N
\right]$ is an optimum of some associated optimization problem.

From the perspective of practical applications, all the goals above are important.
The main interest in this article is to establish the conditions to ensure predictability. To realize the state of the controller, the filter, and the agents, we consider an augmented state space $\mathbb{X} \subset \R^d$. To do this, we use the framework of IRFs, by expressing the issue of predictability due to ergodicity. Sufficient conditions for the existence of a unique, attractive
invariant measure (ergodicity) can be given in terms of  ``{\em average contractivity}''. This key notion can be traced back to \cite{Elton1987,Barnsley1988(2),Barnsley1989}.
The following result is the main idea that leads to the analysis of predictability, fairness, and optimality. One can incorporate \cite[Theorem 2]{Fioravanti2019} with
\cite[Corollary 1]{Elton1987}, to conclude that for all (deterministic) initial conditions $x_0 \in\mathbb X$ and continuous $f:\mathbb X\to \R$, the limit
\begin{equation}\label{eq:ergodicprop}
\lim_{k \to\infty} \frac{1}{k+1} \sum_{j=0}^k  f(X(j)) = \mathbb{E}_\nu(f)
\end{equation}
exists almost surely ($\probability_{x_0}$)  and does not depend on $x_0\in \mathbb X$.  To know more in this direction we refer to \cite{Elton1987,Barnsley1988(2),Barnsley1989,Barnsley2013,Stenflo2001(1),Szarek2003(1),Steinsaltz1999, Walkden2007,Barany2015} and surveys \cite{Diaconis1999,Iosifescu2009,Stenflo2012(s)}.

\section{An Interconnection Result for Two Ensemble Systems}\label{sec:inter-two-ens}
Our first result on the statistical stability of a two-sided market, as modelled by the closed-loop system in Figure \ref{fig:system_ii}, assumes that the agents' states evolve as a linear time-invariant dynamical system, as follows. Let us consider a two-sided market, where one side is modelled by $N_1$ participants
\begin{equation}\label{eq:S1}
{\mathcal S}^1_i ~:~ \left\{ \begin{array}{rcl}
x^{1}_i(k+1) & = & A^{1}_i x^{1}_i(k) + b^{1}_i, \vspace{0.1cm} \\
y^{1}_i(k) & = & {c^{1}_i}^T x^{1}_i(k) + d^{1}_i,
\end{array}\right.
\end{equation}
$i=1,2,\dots, N_1$, responding to the control signal $\pi^1$ produced by
\begin{equation} \label{eq:C1}
{\mathcal C^1} ~:~ \left\{ \begin{array}{rcl}
x^1_c(k+1) & = & A^1_c x_c^1(k) + B^1_c e^1(k), \vspace{0.1cm} \\
\pi^1(k) & = & C^1_c x_c^1(k) + D^1_c e^1(k),
\end{array} \right.
\end{equation}
where the error signal $e^1$ depends on the filtered state of the other side of the market, as follows:
\begin{equation}\label{eq:e1}
e^1(k)  =  u^1(k) - \hat y^2(k).
\end{equation}
The first side has its state filtered too:
\begin{equation}\label{eq:F1}
{\mathcal F^1} ~:~ \left\{ \begin{array}{rcl}
x_f^1(k+1) & = & A^1_f x_f^1(k) + B^1_f y^1(k),\vspace{0.1cm}\\
\hat y^1(k) & = & C^1_f x_f^1(k).
\end{array} \right.
\end{equation}
The other side of the market is modelled by another population of participants
\begin{equation}\label{eq:S2}
{\mathcal S}^2_i ~:~ \left\{ \begin{array}{rcl}
x^{2}_i(k+1) & = & A^{2}_i x^{2}_i(k) + b^{2}_i, \vspace{0.1cm} \\
y^{2}_i(k) & = & {c^{2}_i}^{T} x^{2}_i(k) + d^{2}_i,
\end{array}\right.
\end{equation}
$i=1,2,\dots, N_2$, evolving according to their dynamics and the control of a (in general) different controller
\begin{equation} \label{eq:C2}
{\mathcal C^2} ~:~ \left\{ \begin{array}{rcl}
x^2_c(k+1) & = & A^2_c x_c^2(k) + B^2_c e^2(k), \vspace{0.1cm} \\
\pi^2(k) & = & C^2_c x_c^2(k) + D^2_c e^2(k),
\end{array} \right.
\end{equation}
where the error signal $e^2$ considers the filtered state of the first side of the two-sided market:
\begin{equation}
\label{e2}
e^2(k)  = u^2(k) + \hat y^1(k).
\end{equation}
Likewise, there is a filter for the second side of the two-sided market:
\begin{equation}\label{eq:F2}
{\mathcal F^2} ~:~ \left\{ \begin{array}{rcl}
x_f^2(k+1)  &= & A^2_f x_f^2(k) + B^2_f y(k),\vspace{0.1cm}\\
\hat y^2(k)  &= & C^2_f x_f^2(k).
\end{array} \right.
\end{equation}
\begin{thm} \label{thm:ITE-NLS-LCF}
Consider the feedback system in Figure \ref{fig:system_ii}, with ${\mathcal C}^1, {\mathcal C}^2$ and ${\mathcal F}^1, {\mathcal F}^2$ given in \eqref{eq:C1}, \eqref{eq:C2}, \eqref{eq:F1}, and \eqref{eq:F2}, respectively. 
Assume that each agent  
$i \in \{1,\cdots,N\}$ has its state $x_i$ with dynamics determined by the
equations given in \eqref{eq:S1} and \eqref{eq:S2}, where $A_i^1, A_i^2$ are Schur matrices and $b_i^1, b_i^2$ and $d_i^1, d_i^2$ are taken from the sets $\{b^1_{ij},b^2_{ij}\} \subset \R^{n^1_i}$ and $\{d^1_{i\ell}, d^2_{i\ell}\}
\subset \R^{n^2_i}$ with probability functions $p^1_{ij}(\cdot), p^2_{ij}(\cdot)$, respectively $p^{'1}_{i\ell}(\cdot)$, $p^{'2}_{i\ell}(\cdot)$, that verify \eqref{eq:problaws} and satisfy a Dini continuity condition.
Also, suppose there exist $\delta^1, \delta^2,  \delta^{'1}, \delta^{'2} > 0$ for which
$p^1_{ij}(\pi) \geq  \delta^1 > 0$,
$p^2_{ij}(\pi) \geq  \delta^2 > 0$,
$p^{'1}_{ij}(\pi) \geq  \delta^{'1} > 0$,
$p^{'2}_{ij}(\pi) \geq  \delta^{'2} > 0$ for all
$(i,j)$ and all $\pi^1 \in \Pi^1, \pi^2 \in \Pi^2$. Then, for any stable linear controllers $\mathcal{C}^1$,  $\mathcal{C}^2$ and any stable linear filters $\mathcal{F}^1, \mathcal{F}^2$, the feedback system converges in distribution to a unique invariant measure.

Proof: Following \cite{Barnsley1988(2)},  consider the augmented state
$$
\xi \defeq \left[
((x^1_i)_{i=1}^{N_1})^T \
((x^2_i)_{i=1}^{N_2})^T \
(x_f^1)^T \
(x_f^2)^T \
(x_c^1)^T \
(x_c^2)^T
\right]^T \in \mathbb{X}_{\mathcal S^1} \times \mathbb{X}_{\mathcal S^2} \times \mathbb{X}_{\mathcal F^1} \times \mathbb{X}_{\mathcal F^2} \times \mathbb{X}_{\mathcal C^1} \times \mathbb{X}_{\mathcal C^2},
$$
described by 
$
\xi(k+1) = {\mathcal W}_\ell(x) \defeq \mathcal{A}\xi(k) + \beta_\ell,
$
where $\beta_\ell$ is built from vectors $b^1_{ij}, b^2_{ij}$, scalars ${d}^1_{ij}, {d}^2_{ij}$ and other signals; and
\begin{align}
\label{eq:mathcalA-mat}
\mathcal{A}:=
\begin{bmatrix}
\hat A^1 & 0 & 0 & 0 & 0 & 0\\
0 & \hat A^2 & 0 & 0 & 0 & 0 \\
B^1_f {\bf 1}^T \hat C^1 & 0 & A_f^1 & 0 & 0 & 0 \\
0 & B^2_f {\bf 1}^T \hat C^2 & 0 & A_f^2 & 0 & 0\\
0 & 0 & 0 & -B_c^1 C_f^2 & A_c^1 & 0 \\
0 & 0 & B_c^2 C_f^1 & 0 & 0 & A_c^2 \\
\end{bmatrix}
\end{align}
where $\mathbf{1}$ is the vector of ones, ${\hat A}^1 \defeq \mathbf{diag}(A^1_i)$, ${\hat C}^1 \defeq
\mathbf{diag}(c^{1T}_i)$, ${\hat A}^2 \defeq \mathbf{diag}(A^2_i)$ and ${\hat C}^2 \defeq
\mathbf{diag}(c^{2T}_i)$. To apply Corollary 2.3 from \cite{Barnsley1988(2)}, we make the following observations. First, each map ${\mathcal W}_\ell$ is chosen with probability $p_\ell(\pi) \geq\prod_{i=1}^N \delta_i > 0$. Consequently, these probabilities are bounded away from zero. Second, since $\mathcal A$ is a lower triangular block matrix, it is easy to notice that the spectrum of $\mathcal A$ (i.e., $\sigma\left(\mathcal A\right)$) can be expressed as follows:
$
\sigma({\mathcal A})
= \sigma({\hat A}^1)\cup\sigma({\hat A}^2) \cup\sigma(A^1_f)\cup \sigma(A^2_f)\cup \sigma(A^1_c)\cup \sigma(A^2_c).     
$
By hypothesis, for all $1\le i\le N$, $A_i^1, A_i^2$ are Schur matrices \citep{Schur1921, Graham1976, Bof2018} and  $A_f^1, A_f^2, A_c^1, A_c^2$ are Schur matrices as well. For any induced matrix norm, there exists $m
\in \N$ sufficiently large such that $\|{\mathcal A}^m\| < 1$. The result then is a direct consequence of \cite{Barnsley1988(2)}. $\blacksquare$
\end{thm}
In our next result, the agents' states are realized by a nonlinear iterated random function.

\begin{thm} \label{thm:intrcon-twoensmbl-nonlinrstate-lincontr&filtr}
Consider the system shown in Figure \ref{fig:system_ii}.
Assume that each agent
$i \in \{1,\cdots,N_1\}$ and  $i \in \{1,\cdots,N_2\}$ has a state realized by the nonlinear iterated random functions
\begin{align}
\label{eq:nonlinear1}
x^1_i(k+1) &= {\mathcal W}^1_{ij}(x^1_i(k)),
\end{align}
\begin{align}
y^1_i(k) &= {\mathcal H}^1_{ij}(x^1_i(k)),
\end{align}
and 
\begin{align}
\label{eq:nonlinear2}
x^2_i(k+1) &= {\mathcal W}^2_{ij}(x^2_i(k)),
\end{align}
\begin{align}
\label{eq:nonlinear4}
y^2_i(k) &= {\mathcal H}^2_{ij}(x^2_i(k)),
\end{align}
where ${\mathcal W}^1_{ij}, {\mathcal H}^1_{ij}, {\mathcal W}^2_{ij}, {\mathcal H}^2_{ij}$ are
globally Lipschitz-continuous functions with
global Lipschitz constants $l^1_{ij}, l^{'1}_{ij}, l^2_{ij}, l^{'2}_{ij}$ respectively.
In addition, we assume that we are dealing with Dini continuous probability functions
$p^1_{ij},p^{'1}_{il}, p^2_{ij},p^{'2}_{il}$ so that \eqref{eq:problaws} are satisfied. We also assume
there exist  scalars $\delta^1, \delta^{'1},\delta^2, \delta^{'2}> 0$ so that
$p^1_{ij}(\pi) \geq  \delta^1 > 0$,
$p^{'1}_{ij}(\pi) \geq  \delta^{'1} > 0$,
$p^2_{ij}(\pi) \geq  \delta^2 > 0$,
$p^{'2}_{ij}(\pi) \geq  \delta^{'2} > 0$ 
for all $(i,j)$.
\noindent
Furthermore, assume that the following contractivity condition holds:
for all $1 \le i \le N_1, 1 \le j \le J_1, 1 \le i \le N_2, 1 \le j \le J_2$, we have $l^1_{ij} < 1, l^2_{ij}<1$.
\noindent
Then, for every linear controller with eigenvalues in the unit circle $\mathcal{C}_1$ and $\mathcal {C}_2$ and every
linear filter with eigenvalues in the unit circle $\mathcal{F}_1$ and $\mathcal{F}_2$ compatible with the feedback structure in Figure \ref{fig:system_ii}, the feedback loop has a unique attractive invariant measure and is ergodic.
\end{thm}
Proof: Similar to the proof of \cite[Theorem 18]{Fioravanti2019}, the assumptions on the
Lipschitz constants and the internal asymptotic stability of
controller and filter permit the application of Theorem 2.1 and Corollary 2.2 of \cite{Barnsley1988(2)}. $\blacksquare$


We now consider ensembles where the agents' actions are limited to a finite set. In this case, the ergodic behaviour follows from the results of  \cite{Werner2005, Werner2004}.

\begin{thm}\label{thm:discrete-range-space}
Consider the system given in Figure \ref{fig:system_ii}.  Assume
that ${\mathbb A}_i$ is finite for all $i$ and that agent $i$ 
has a state realized by the nonlinear stochastic
difference equations (\ref{eq:nonlinear1}--\ref{eq:nonlinear4}). 
\noindent
In addition, as before, assume Dini continuous probability functions $p^1_{ij},p^{1'}_{il}, p^2_{ij},p^{2'}_{il}$ so that
\eqref{eq:problaws} hold. Assume
furthermore that there are scalars $\delta^1, \delta^{1'}, \delta^2, \delta^{2'} > 0$ such that
$p^1_{ij}(\pi) \geq \delta^1 > 0$, $p^{1'}_{ij}(\pi) \geq \delta^{1'} > 0$, $p^2_{ij}(\pi) \geq \delta^2 > 0$, $p^{2'}_{ij}(\pi) \geq \delta^{2'} > 0$  for all
$(i,j)$ and all $\pi$. Then, for all linear controller $\mathcal{C}$ with eigenvalues inside the unit circle and
every linear filter with eigenvalues inside the unity circle, the following holds:

If the graph $\mathcal G=(\mathbb{X}_S, E)$ is strongly connected 
then there exists an
invariant measure feedback system. If, in addition, the adjacency matrix of
the graph is primitive, then the invariant measure is uniquely attractive.
\end{thm}

Proof: This is a consequence of \cite{Werner2004} and the observation
that the necessary contractivity properties follow from the internal
asymptotic stability of the controller and filter. $\blacksquare$


\section{A Large-scale Interconnection for Ensemble Systems}\label{sec:inter-lrgscl-ensmbl}
As hinted at in Section~\ref{app:multi-sided}, many markets traditionally thought of as two-sided are, in fact, multi-sided platforms \citep{Weyl2010,Hagiu2015,Evans2016}.
In the example of Uber, for instance, one can also consider ``fleet partners'' (who are intermediaries for car manufacturers and car leasing providers) and taxi operators. 
Such interconnections can get rather complicated, and we would like to capture these interconnections in a matrix, similar to the adjacency matrix of a graph. 

Thus, building on the techniques of Section \ref{sec:inter-two-ens}, consider the interconnection of a multi-sided platform described as follows:
\begin{equation}\label{eq:mult-side}
e^p = u^p - \sum^M_{q=1}H_{pq}\hat{y}^q,
\end{equation}
where $e^p$, for $p = 1,2,\dots,M$, is the input to 
\begin{equation} \label{eq:Cp}
{\mathcal C^p} ~:~ \left\{ \begin{array}{rcl}
x^p_c(k+1) & = & A^p_c x_c^p(k) + B^p_c e^p(k), \vspace{0.1cm} \\
\pi^p(k) & = & C^p_c x_c^p(k) + D^p_c e^p(k),
\end{array} \right.
\end{equation}
$u^p$ are external inputs, $H_{pq}$ is a matrix with real and constant entries, and $\hat{y}^q$, for $q = 1,2,\dots,M$, is the output of
\begin{equation}\label{eq:Fq}
{\mathcal F^q} ~:~ \left\{ \begin{array}{rcl}
x_f^q(k+1) & = & A^q_f x_f^q(k) + B^q_f y^q(k),\vspace{0.1cm}\\
\hat y^q(k) & = & C^q_f x_f^q(k).
\end{array} \right.
\end{equation}
We also have, for $i=1,2,\dots,N_q$,
\begin{equation}\label{eq:Sq}
{\mathcal S}^q_i ~:~ \left\{ \begin{array}{rcl}
x^{q}_i(k+1) & = & A^{q}_i x^{q}_i(k) + b^{q}_i, \vspace{0.1cm} \\
y^{q}_i(k) & = & {c^{q}_i}^T x^{q}_i(k) + d^{q}_i,
\end{array}\right.
\end{equation}
where $N_q$ is the total number of systems in the $q$th ensemble denoted by ${\mathcal S}^q$. Finally, $y^{q}(k) = y^{q}_1(k) + y^{q}_2(k) + \dots + y^{q}_{N_q}(k)$.

\begin{figure}[ht]
\centering
\includegraphics[width=\columnwidth]{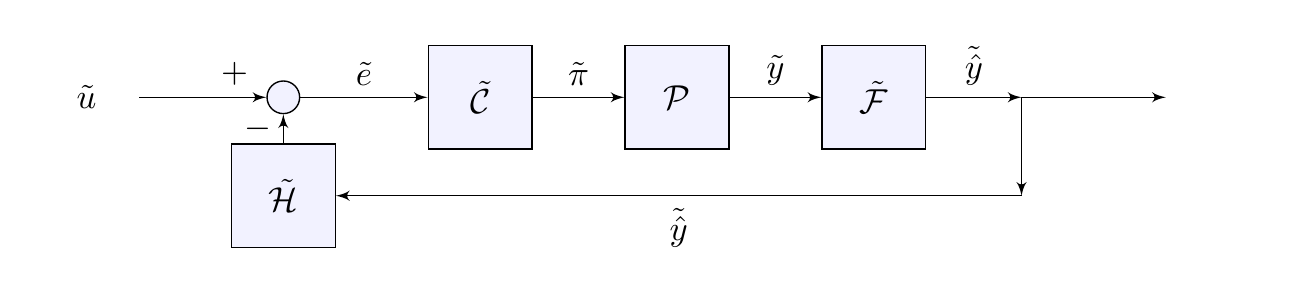}
\caption{A feedback model for the large-scale interconnection.}
\label{fig:largescale}
\end{figure}

Now, by writing
\begin{equation*}
\tilde{e}:=\begin{bmatrix}e^1 \\ e^2 \\ \vdots \\ e^M\end{bmatrix}, \
\tilde{u}:=\begin{bmatrix}u^1 \\ u^2 \\ \vdots \\ u^M\end{bmatrix}, \
\tilde{\pi}:=\begin{bmatrix}\pi^1 \\ \pi^2 \\ \vdots \\ \pi^M\end{bmatrix}, \
\tilde{y}:=\begin{bmatrix}y^1 \\ y^2 \\ \vdots \\ y^M\end{bmatrix} \text{ and }
\tilde{\hat{y}}:=\begin{bmatrix}\hat{y}^1 \\ \hat{y}^2 \\ \vdots \\ \hat{y}^M\end{bmatrix}, 
\end{equation*}
the interconnection description may be expressed more compactly as
\begin{align}\label{eq:int-descrptn}
\Tilde{e}=\Tilde{u}-{\Tilde{\mathcal H}} \Tilde{\hat{y}},
\end{align}
where $\Tilde{\mathcal H}$ is a matrix with block entries $H_{pq}$. Let $\Tilde{\mathcal C}:=\text{diag}(\mathcal{C}^1,\dots,\mathcal{C}^M)$ such that $\Tilde{\pi}= \Tilde{\mathcal C} \Tilde{e}$, and let $\Tilde{\mathcal F}:=\text{diag}(\mathcal{F}^1,\dots,\mathcal{F}^M)$ such that 
$\Tilde{\hat{y}}= \Tilde{\mathcal F} \Tilde{y}$. This set-up is depicted in Figure \ref{fig:largescale}. Then, similar to Theorem \ref{thm:ITE-NLS-LCF}, we have the following result.

\begin{thm} \label{thm:lrg-scale-intr-cnnctn}
Consider the feedback system described in the preceding two paragraphs and depicted in Figure \ref{fig:largescale}. For each ensemble ${\mathcal S}^q$, $q = 1,2,\dots, M$, assume that each $i$th system in the ensemble has its state $x^q_i$ with dynamics determined by the affine stochastic difference equations given by \eqref{eq:Sq}, where $A_i^q$ are Schur matrices, and $b_i^q$ and $d_i^q$ are chosen, at each time step, from sets according to Dini continuous probability functions in the manner of Theorem \ref{thm:ITE-NLS-LCF}, and that these probability functions are bounded below by scalars strictly more significant than $0$. Then, for any stable linear controllers  $\mathcal{C}^1,\dots,\mathcal{C}^M$ and any stable linear filters $\mathcal{F}^1,\dots,\mathcal{F}^M$ compatible with the system structure, the feedback loop converges in distribution to a unique invariant measure.
\end{thm}
Proof: In the spirit of Theorem \ref{thm:ITE-NLS-LCF}, by defining an augmented state
\begin{equation*}
\xi \defeq \begin{bmatrix}
\Tilde{x}\\
\Tilde{x}_f\\
\Tilde{x}_c
\end{bmatrix}
\end{equation*}
whose dynamics are realized by the following equation:
\begin{align*}
\xi(k+1) = {\mathcal W}_\ell(x) \defeq \mathcal{A}\xi(k) + \beta_\ell,       
\end{align*}
where
\begin{equation*}
\Tilde{x}_f \defeq \begin{bmatrix}
x^1_f \\ x^2_f \\ \vdots \\x^M_f
\end{bmatrix}, \
\Tilde{x}_c \defeq \begin{bmatrix}
x^1_c \\ x^2_c \\ \vdots \\x^M_c
\end{bmatrix}, \
\Tilde{x} \defeq \begin{bmatrix}
x^1 \\ x^2 \\ \vdots \\x^M
\end{bmatrix}, \
x^1 \defeq \begin{bmatrix}
x^1_1 \\ x^1_2 \\ \vdots \\x^1_{N_1}
\end{bmatrix}, \
x^2 \defeq \begin{bmatrix}
x^2_1 \\ x^2_2 \\ \vdots \\x^2_{N_2}
\end{bmatrix}, \ \dots, \
x^M \defeq \begin{bmatrix}
x^M_1 \\ x^M_2 \\ \vdots \\x^M_{N_M}
\end{bmatrix},
\end{equation*}
we obtain
\begin{equation}\label{mathcalA-mat}
\mathcal{A}:=\begin{bmatrix}
\hat{A}&0&0\\
B_f\hat{\mathbf{1}}\hat{C}& A_f&0\\
0& -B_c \Tilde{\mathcal H} C_f& A_c
\end{bmatrix}.
\end{equation}
In \eqref{mathcalA-mat}, $\mathbf{1}$ is the vector of ones, $\hat{\mathbf{1}} := \text{diag}(\mathbf{1}^T)$ and
\begin{gather*}
\hat{A} := \text{diag}(A^1_1, A^1_2, \dots, A^1_{N_1}, A^2_1, A^2_2, \dots, A^2_{N_2}, \dots, A^M_1, A^M_2, \dots, A^M_{N_M}),\\
\hat{C} := \text{diag}(c^{1T}_1, c^{1T}_2, \dots, c^{1T}_{N_1}, c^{2T}_1, c^{2T}_2, \dots, c^{2T}_{N_2}, \dots, c^{MT}_1, c^{MT}_2, \dots, c^{MT}_{N_M}),\\
A_f := \text{diag}(A^1_f, A^2_f, \dots, A^M_f), \
B_f  := \text{diag}(B^1_f, B^2_f, \dots, B^M_f), \
C_f := \text{diag}(C^1_f, C^2_f, \dots, C^M_f),\\
A_c := \text{diag}(A^1_c, A^2_c, \dots, A^M_c), \
B_c := \text{diag}(B^1_c, B^2_c, \dots, B^M_c).
\end{gather*}
The proof then follows in a manner similar to the proof of Theorem \ref{thm:ITE-NLS-LCF}. $\blacksquare$

\section{Conclusions and Further Work}\label{sec:conclsn}
In feedback control systems, a demanding and emerging area for further study is the control of ensembles of agents. In practice, one such example is an online labour platform \citep{mohlmann2020algorithmic}.
There are two main differences between the control of ensemble problems from the classical control problems. First, although the ensembles generally are too large to allow for a microscopic approach, they are not sufficiently large to allow for a meaningful fluid (mean-field) approximation.
Second, the regulation problem concerns the ensemble and the individual agents; a certain quality of service should be provided to each agent.
We have formulated this problem as an iterated random function to design an ergodic control which is the key to delivering the expected quality of service to the agents across the network.


\bibliographystyle{ieeetran}
\bibliography{ref,mdps,twosided}

\begin{thebibliography}{78}
\providecommand{\natexlab}[1]{#1}
\providecommand{\url}[1]{#1}
\csname url@samestyle\endcsname
\providecommand{\newblock}{\relax}
\providecommand{\bibinfo}[2]{#2}
\providecommand{\BIBentrySTDinterwordspacing}{\spaceskip=0pt\relax}
\providecommand{\BIBentryALTinterwordstretchfactor}{4}
\providecommand{\BIBentryALTinterwordspacing}{\spaceskip=\fontdimen2\font plus
\BIBentryALTinterwordstretchfactor\fontdimen3\font minus
  \fontdimen4\font\relax}
\providecommand{\BIBforeignlanguage}[2]{{%
\expandafter\ifx\csname l@#1\endcsname\relax
\typeout{** WARNING: IEEEtranN.bst: No hyphenation pattern has been}%
\typeout{** loaded for the language `#1'. Using the pattern for}%
\typeout{** the default language instead.}%
\else
\language=\csname l@#1\endcsname
\fi
#2}}
\providecommand{\BIBdecl}{\relax}
\BIBdecl

\bibitem[Rochet and Tirole(2003)]{Rochet2003}
J.-C. Rochet and J.~Tirole, ``{Platform Competition in Two-Sided Markets},''
  \emph{Journal of the European Economic Association}, vol.~1, no.~4, pp.
  990--1029, 06 2003.

\bibitem[Rochet and Tirole(2004{\natexlab{a}})]{Rochet2004(1)}
------, ``Defining two-sided markets,'' Citeseer, Tech. Rep., 2004.

\bibitem[Rochet and Tirole(2004{\natexlab{b}})]{Rochet2004(2)}
------, ``Two-sided markets: an overview,'' \emph{Institut d’Economie
  Industrielle working paper}, 2004.

\bibitem[Rochet and Tirole(2006)]{Rochet2006}
------, ``Two-sided markets: a progress report,'' \emph{The RAND journal of
  economics}, vol.~37, no.~3, pp. 645--667, 2006.

\bibitem[Weyl(2010)]{Weyl2010}
\BIBentryALTinterwordspacing
E.~G. Weyl, ``A price theory of multi-sided platforms,'' \emph{American
  Economic Review}, vol. 100, no.~4, pp. 1642--72, September 2010. [Online].
  Available: \url{https://www.aeaweb.org/articles?id=10.1257/aer.100.4.1642}
\BIBentrySTDinterwordspacing

\bibitem[Hagiu and Wright(2015)]{Hagiu2015}
A.~Hagiu and J.~Wright, ``Multi-sided platforms,'' \emph{International Journal
  of Industrial Organization}, vol.~43, pp. 162--174, 2015.

\bibitem[Evans and Schmalensee(2016)]{Evans2016}
D.~S. Evans and R.~Schmalensee, \emph{Matchmakers: The new economics of
  multisided platforms}.\hskip 1em plus 0.5em minus 0.4em\relax Harvard
  Business Review Press, 2016.

\bibitem[Parker et~al.(2016{\natexlab{a}})Parker, Van~Alstyne, and
  Jiang]{Parker2016(1)}
G.~Parker, M.~W. Van~Alstyne, and X.~Jiang, ``Platform ecosystems: How
  developers invert the firm,'' \emph{MIS Quarterly}, vol.~41, no.~1, pp.
  255--266, 2016.

\bibitem[Ashlagi et~al.(2020)Ashlagi, Braverman, Kanoria, and
  Shi]{ashlagi2020clearing}
I.~Ashlagi, M.~Braverman, Y.~Kanoria, and P.~Shi, ``Clearing matching markets
  efficiently: informative signals and match recommendations,''
  \emph{Management Science}, vol.~66, no.~5, pp. 2163--2193, 2020.

\bibitem[Arnosti et~al.(2021)Arnosti, Johari, and Kanoria]{Arnosti2021}
N.~Arnosti, R.~Johari, and Y.~Kanoria, ``Managing congestion in matching
  markets,'' \emph{Manufacturing \& Service Operations Management}, vol.~23,
  no.~3, pp. 620--636, 2021.

\bibitem[Benjaafar and Hu(2021)]{Benjaafar2021}
\BIBentryALTinterwordspacing
S.~Benjaafar and M.~Hu, ``Introduction to the special issue on sharing economy
  and innovative marketplaces,'' \emph{Manufacturing \& Service Operations
  Management}, vol.~23, no.~3, pp. 549--552, 2021. [Online]. Available:
  \url{https://doi.org/10.1287/msom.2021.0998}
\BIBentrySTDinterwordspacing

\bibitem[Liu et~al.(2021)Liu, Zhang, and Zhang]{Liu2021}
\BIBentryALTinterwordspacing
Z.~Liu, D.~J. Zhang, and F.~Zhang, ``Information sharing on retail platforms,''
  \emph{Manufacturing \& Service Operations Management}, vol.~23, no.~3, pp.
  606--619, 2021. [Online]. Available:
  \url{https://pubsonline.informs.org/doi/abs/10.1287/msom.2020.0915}
\BIBentrySTDinterwordspacing

\bibitem[Kanoria et~al.(2018)Kanoria, Saban, and
  Sethuraman]{kanoria2018convergence}
Y.~Kanoria, D.~Saban, and J.~Sethuraman, ``Convergence of the core in
  assignment markets,'' \emph{Operations Research}, vol.~66, no.~3, pp.
  620--636, 2018.

\bibitem[M{\"o}hlmann et~al.(2021)M{\"o}hlmann, Zalmanson, Henfridsson, and
  Gregory]{mohlmann2020algorithmic}
M.~M{\"o}hlmann, L.~Zalmanson, O.~Henfridsson, and R.~W. Gregory, ``Algorithmic
  management of work on online labor platforms: When matching meets control,''
  \emph{MIS Quarterly}, vol.~45, 2021.

\bibitem[Lundy et~al.(2019)Lundy, Wei, Fu, Kominers, and
  Leyton-Brown]{Lundy2019}
\BIBentryALTinterwordspacing
T.~Lundy, A.~Wei, H.~Fu, S.~D. Kominers, and K.~Leyton-Brown, ``Allocation for
  social good: Auditing mechanisms for utility maximization,'' in
  \emph{Proceedings of the 2019 ACM Conference on Economics and Computation},
  ser. EC '19.\hskip 1em plus 0.5em minus 0.4em\relax New York, NY, USA:
  Association for Computing Machinery, 2019, p. 785–803. [Online]. Available:
  \url{https://doi.org/10.1145/3328526.3329623}
\BIBentrySTDinterwordspacing

\bibitem[Fioravanti et~al.(2019)Fioravanti, Marecek, Shorten, Souza, and
  Wirth]{Fioravanti2019}
A.~R. Fioravanti, J.~Marecek, R.~N. Shorten, M.~Souza, and F.~R. Wirth, ``On
  the ergodic control of ensembles,'' \emph{Automatica}, vol. 108, p. 108483,
  2019.

\bibitem[Bateni et~al.(2016)Bateni, Chen, Ciocan, and Mirrokni]{Bateni2016}
\BIBentryALTinterwordspacing
M.~H. Bateni, Y.~Chen, D.~F. Ciocan, and V.~Mirrokni, ``Fair resource
  allocation in a volatile marketplace,'' in \emph{Proceedings of the 2016 ACM
  Conference on Economics and Computation}, ser. EC '16.\hskip 1em plus 0.5em
  minus 0.4em\relax New York, NY, USA: Association for Computing Machinery,
  2016, p. 819. [Online]. Available:
  \url{https://doi.org/10.1145/2940716.2940763}
\BIBentrySTDinterwordspacing

\bibitem[Lobel(2020)]{Lobel2020}
I.~Lobel, ``Revenue management and the rise of the algorithmic economy,''
  \emph{Management Science}, 2020.

\bibitem[Chen(2016)]{Chen2016}
\BIBentryALTinterwordspacing
M.~K. Chen, ``Dynamic pricing in a labor market: Surge pricing and flexible
  work on the uber platform,'' in \emph{Proceedings of the 2016 ACM Conference
  on Economics and Computation}, ser. EC '16.\hskip 1em plus 0.5em minus
  0.4em\relax New York, NY, USA: Association for Computing Machinery, 2016, p.
  455. [Online]. Available: \url{https://doi.org/10.1145/2940716.2940798}
\BIBentrySTDinterwordspacing

\bibitem[Castillo et~al.(2017)Castillo, Knoepfle, and Weyl]{Castillo2017}
J.~C. Castillo, D.~Knoepfle, and G.~Weyl, ``Surge pricing solves the wild goose
  chase,'' in \emph{Proceedings of the 2017 ACM Conference on Economics and
  Computation}, 2017, pp. 241--242.

\bibitem[Cachon et~al.(2017)Cachon, Daniels, and Lobel]{Cachon2017}
G.~P. Cachon, K.~M. Daniels, and R.~Lobel, ``The role of surge pricing on a
  service platform with self-scheduling capacity,'' \emph{Manufacturing \&
  Service Operations Management}, vol.~19, no.~3, pp. 368--384, 2017.

\bibitem[Castillo(2020)]{Castillo2020}
J.~C. Castillo, ``Who benefits from surge pricing?'' \emph{Available at SSRN
  3245533}, 2020.

\bibitem[Garg and Nazerzadeh(2020)]{Garg2020}
\BIBentryALTinterwordspacing
N.~Garg and H.~Nazerzadeh, ``Driver surge pricing,'' in \emph{Proceedings of
  the 21st ACM Conference on Economics and Computation}, ser. EC '20.\hskip 1em
  plus 0.5em minus 0.4em\relax New York, NY, USA: Association for Computing
  Machinery, 2020, p. 501. [Online]. Available:
  \url{https://doi.org/10.1145/3391403.3399476}
\BIBentrySTDinterwordspacing

\bibitem[Simonetto et~al.(2019)Simonetto, Monteil, and Gambella]{Simonetto2019}
A.~Simonetto, J.~Monteil, and C.~Gambella, ``Real-time city-scale ridesharing
  via linear assignment problems,'' \emph{Transportation Research Part C:
  Emerging Technologies}, vol. 101, pp. 208 -- 232, 2019.

\bibitem[Araman et~al.(2019)Araman, Calmon, and Fridgeirsdottir]{Araman2019}
V.~F. Araman, A.~Calmon, and K.~Fridgeirsdottir, ``Pricing and job allocation
  in online labor platforms,'' \emph{INSEAD Working Paper No. 2019/32/TOM},
  2019.

\bibitem[Aouad and Sarita{\c{c}}(2020)]{Aouad2020}
A.~Aouad and {\"O}.~Sarita{\c{c}}, ``Dynamic stochastic matching under limited
  time,'' in \emph{Proceedings of the 21st ACM Conference on Economics and
  Computation}, 2020, pp. 789--790.

\bibitem[{\"O}zkan(2020)]{Ozkan2020}
E.~{\"O}zkan, ``Joint pricing and matching in ride-sharing systems,''
  \emph{European Journal of Operational Research}, vol. 287, no.~3, pp.
  1149--1160, 2020.

\bibitem[Yan et~al.(2020)Yan, Zhu, Korolko, and Woodard]{Yan2020}
C.~Yan, H.~Zhu, N.~Korolko, and D.~Woodard, ``Dynamic pricing and matching in
  ride-hailing platforms,'' \emph{Naval Research Logistics (NRL)}, vol.~67,
  no.~8, pp. 705--724, 2020.

\bibitem[Parker et~al.(2016{\natexlab{b}})Parker, Van~Alstyne, and
  Choudary]{Parker2016(2)}
G.~G. Parker, M.~W. Van~Alstyne, and S.~P. Choudary, \emph{Platform revolution:
  How networked markets are transforming the economy and how to make them work
  for you}.\hskip 1em plus 0.5em minus 0.4em\relax WW Norton \& Company, 2016.

\bibitem[Ma et~al.(2020)Ma, Fang, and Parkes]{Ma2020}
H.~Ma, F.~Fang, and D.~C. Parkes, ``Spatio-temporal pricing for ridesharing
  platforms,'' \emph{ACM SIGecom Exchanges}, vol.~18, no.~2, pp. 53--57, 2020.

\bibitem[Cook et~al.(2018)Cook, Diamond, Hall, List, and Oyer]{Cook2018}
C.~Cook, R.~Diamond, J.~Hall, J.~A. List, and P.~Oyer, ``The gender earnings
  gap in the gig economy: Evidence from over a million rideshare drivers,''
  National Bureau of Economic Research, Tech. Rep., 2018.

\bibitem[S{\"u}hr et~al.(2019)S{\"u}hr, Biega, Zehlike, Gummadi, and
  Chakraborty]{Suhr2019}
T.~S{\"u}hr, A.~J. Biega, M.~Zehlike, K.~P. Gummadi, and A.~Chakraborty,
  ``Two-sided fairness for repeated matchings in two-sided markets: A case
  study of a ride-hailing platform,'' in \emph{Proceedings of the 25th ACM
  SIGKDD International Conference on Knowledge Discovery \& Data Mining}, 2019,
  pp. 3082--3092.

\bibitem[Cohen et~al.(2019)Cohen, Elmachtoub, and Lei]{Cohen2019}
M.~Cohen, A.~N. Elmachtoub, and X.~Lei, ``Price discrimination with fairness
  constraints,'' \emph{Available at SSRN 3459289}, 2019.

\bibitem[Jung et~al.(2020)Jung, Kannan, Lee, Pai, Roth, and Vohra]{Jung2020}
\BIBentryALTinterwordspacing
C.~Jung, S.~Kannan, C.~Lee, M.~Pai, A.~Roth, and R.~Vohra, ``Fair prediction
  with endogenous behavior,'' in \emph{Proceedings of the 21st ACM Conference
  on Economics and Computation}, ser. EC '20.\hskip 1em plus 0.5em minus
  0.4em\relax New York, NY, USA: Association for Computing Machinery, 2020, p.
  677–678. [Online]. Available: \url{https://doi.org/10.1145/3391403.3399473}
\BIBentrySTDinterwordspacing

\bibitem[Freeman et~al.(2020)Freeman, Shah, and Vaish]{Freeman2020}
\BIBentryALTinterwordspacing
R.~Freeman, N.~Shah, and R.~Vaish, ``Best of both worlds: Ex-ante and ex-post
  fairness in resource allocation,'' in \emph{Proceedings of the 21st ACM
  Conference on Economics and Computation}, ser. EC '20.\hskip 1em plus 0.5em
  minus 0.4em\relax New York, NY, USA: Association for Computing Machinery,
  2020, p. 21–22. [Online]. Available:
  \url{https://doi.org/10.1145/3391403.3399537}
\BIBentrySTDinterwordspacing

\bibitem[Chouldechova and Roth(2020)]{Chouldechova2018}
A.~Chouldechova and A.~Roth, ``A snapshot of the frontiers of fairness in
  machine learning,'' \emph{Commun. ACM}, vol.~63, no.~5, p. 82–89, Apr.
  2020.

\bibitem[Mouzannar et~al.(2019)Mouzannar, Ohannessian, and
  Srebro]{Mouzannar2019}
H.~Mouzannar, M.~I. Ohannessian, and N.~Srebro, ``From fair decision making to
  social equality,'' in \emph{Proceedings of the Conference on Fairness,
  Accountability, and Transparency}, 2019, pp. 359--368.

\bibitem[McArthur et~al.(2007)McArthur, Davidson, Catterson, Dimeas,
  Hatziargyriou, Ponci, and Funabashi]{Agents07}
S.~McArthur, E.~Davidson, V.~Catterson, A.~Dimeas, N.~Hatziargyriou, F.~Ponci,
  and T.~Funabashi, ``Multi-agent systems for power engineering applications
  --- {P}art {I}: {C}oncepts, approaches, and technical challenges,''
  \emph{IEEE Transactions on Power Systems}, vol.~22, no.~4, pp. 1743--1752,
  Nov. 2007.

\bibitem[Blondel et~al.(2005)Blondel, Hendrickx, Olshevsky, and
  Tsitsiklis]{Blondel05}
V.~Blondel, J.~Hendrickx, A.~Olshevsky, and J.~Tsitsiklis, ``Convergence in
  multiagent coordination, consensus, and flocking,'' in \emph{Proceedings of
  the 44th IEEE Conference on Decision and Control, and the European Control
  Conference 2005}, Seville, Spain, Dec. 2005, pp. 2996--3000.

\bibitem[Nedi\'{c} and Ozdaglar(2009)]{Nedic09}
A.~Nedi\'{c} and A.~Ozdaglar, ``Distributed subgradient methods for multi-agent
  optimization,'' \emph{IEEE Transactions on Automatic Control}, vol.~54,
  no.~1, pp. 48--61, Jan. 2009.

\bibitem[Mathew and Mezi\'{c}(2011)]{Mezic11}
G.~Mathew and I.~Mezi\'{c}, ``Metrics for ergodicity and design of ergodic
  dynamics for multi-agent systems,'' \emph{Physica D}, vol. 240, pp. 432--442,
  2011, physica D.

\bibitem[Elton(1987)]{Elton1987}
J.~H. Elton, ``An ergodic theorem for iterated maps,'' \emph{Ergodic Theory and
  Dynamical Systems}, vol.~7, no.~04, pp. 481--488, 1987.

\bibitem[Barnsley and Elton(1988)]{Barnsley1988(1)}
M.~F. Barnsley and J.~H. Elton, ``A new class of markov processes for image
  encoding,'' \emph{Advances in applied probability}, vol.~20, no.~1, pp.
  14--32, 1988.

\bibitem[Barnsley et~al.(1989)Barnsley, Elton, and Hardin]{Barnsley1989}
M.~F. Barnsley, J.~H. Elton, and D.~P. Hardin, ``Recurrent iterated function
  systems,'' \emph{Constructive approximation}, vol.~5, no.~1, pp. 3--31, 1989.

\bibitem[Barnsley et~al.(1988)Barnsley, Demko, Elton, and
  Geronimo]{Barnsley1988(2)}
M.~F. Barnsley, S.~G. Demko, J.~H. Elton, and J.~S. Geronimo, ``Invariant
  measures for markov processes arising from iterated function systems with
  place-dependent probabilities,'' in \emph{Annales de l'IHP Probabilit{\'e}s
  et statistiques}, vol.~24, 1988, pp. 367--394.

\bibitem[Barnsley(2013)]{Barnsley2013}
M.~Barnsley, \emph{Fractals Everywhere}, ser. Dover books on mathematics.\hskip
  1em plus 0.5em minus 0.4em\relax Dover Publications, 2013.

\bibitem[Stenflo(2001{\natexlab{a}})]{Stenflo2001(1)}
{\"O}.~Stenflo, ``Markov chains in random environments and random iterated
  function systems,'' \emph{Transactions of the American Mathematical Society},
  vol. 353, no.~9, pp. 3547--3562, 2001.

\bibitem[Szarek(2003{\natexlab{a}})]{Szarek2003(1)}
T.~Szarek, ``Invariant measures for markov operators with application to
  function systems,'' \emph{Studia Mathematica}, vol. 154, pp. 207--222, 01
  2003.

\bibitem[Steinsaltz(1999)]{Steinsaltz1999}
D.~Steinsaltz, ``Locally contractive iterated function systems,'' \emph{Ann.
  Probab.}, vol.~27, no.~4, pp. 1952--1979, 10 1999.

\bibitem[Walkden(2007)]{Walkden2007}
\BIBentryALTinterwordspacing
C.~P. Walkden, ``Invariance principles for interated maps that contract on
  average,'' \emph{Transactions of the American Mathematical Society}, vol.
  359, no.~3, pp. 1081--1097, 2007. [Online]. Available:
  \url{http://www.jstor.org/stable/20161616}
\BIBentrySTDinterwordspacing

\bibitem[B{\'a}r{\'a}ny(2015)]{Barany2015}
B.~B{\'a}r{\'a}ny, ``On iterated function systems with place-dependent
  probabilities,'' \emph{Proc. Amer. Math. Soc.}, vol. 143, pp. 419--432, 2015.

\bibitem[Diaconis and Freedman(1999)]{Diaconis1999}
P.~Diaconis and D.~Freedman, ``Iterated random functions,'' \emph{SIAM Review},
  vol.~41, no.~1, pp. 45--76, 1999.

\bibitem[Iosifescu(2009)]{Iosifescu2009}
M.~Iosifescu, ``Iterated function systems: A critical survey,'' \emph{Math.
  Reports}, vol.~11, no.~3, pp. 181--229, 2009.

\bibitem[Stenflo(2012)]{Stenflo2012(s)}
{\"O}.~Stenflo, ``A survey of average contractive iterated function systems,''
  \emph{Journal of Difference Equations and Applications}, vol.~18, no.~8, pp.
  1355--1380, 2012.

\bibitem[Schlote et~al.(2013)Schlote, H\"ausler, Hecker, Bergmann, Crisostomi,
  Radusch, and Shorten]{Arieh13}
A.~Schlote, F.~H\"ausler, T.~Hecker, A.~Bergmann, E.~Crisostomi, I.~Radusch,
  and R.~Shorten, ``Cooperative regulation and trading of emissions using
  plug-in hybrid vehicles,'' \emph{IEEE Transactions on Intelligent
  Transportation Systems}, vol.~14, no.~4, pp. 1572--1585, 2013.

\bibitem[Schlote et~al.(2014)Schlote, King, Crisostomi, and Shorten]{Arieh14}
A.~Schlote, C.~King, E.~Crisostomi, and R.~Shorten, ``Delay-tolerant stochastic
  algorithms for parking space assignment,'' \emph{IEEE Transactions on
  Intelligent Transportation Systems}, vol.~15, no.~5, pp. 1922--1935, Oct
  2014.

\bibitem[Marecek et~al.(2015)Marecek, Shorten, and Yu]{marevcek2015signaling}
J.~Marecek, R.~Shorten, and J.~Y. Yu, ``Signaling and obfuscation for
  congestion control,'' \emph{International Journal of Control}, vol.~88,
  no.~10, pp. 2086--2096, 2015.

\bibitem[Stenflo(1998)]{Stenflo1998}
{\"O}.~Stenflo, \emph{Ergodic theorems for time-dependent random iteration of
  functions}.\hskip 1em plus 0.5em minus 0.4em\relax University of Ume{\aa},
  Department of Mathematics, 1998.

\bibitem[Stenflo(1999)]{Stenflo1999}
------, ``Ergodic theorems for iterated function systems controlled by
  stochastic sequences,'' \emph{PhD thesis}, 1999.

\bibitem[Stenflo(2001{\natexlab{b}})]{Stenflo2001(2)}
------, ``A note on a theorem of karlin,'' \emph{Statistics \& probability
  letters}, vol.~54, no.~2, pp. 183--187, 2001.

\bibitem[Szarek(1997)]{Szarek1997}
T.~Szarek, ``Iterated function systems depending on previous transformation,''
  \emph{Acta Mathematica}, 03 1997.

\bibitem[Szarek(1999)]{Szarek1999}
------, ``Generic properties of continuous iterated function systems,''
  \emph{Bulletin of the Polish Academy of Sciences. Mathematics}, vol.~47,
  no.~1, pp. 77--89, 1999.

\bibitem[Szarek(2000{\natexlab{a}})]{Szarek2000(1)}
------, ``Invariant measures for iterated function systems,'' in \emph{Annales
  Polonici Mathematici 75(1)}, 01 2000.

\bibitem[Szarek(2000{\natexlab{b}})]{Szarek2000(2)}
------, ``The stability of markov operators on polish spaces,'' \emph{Studia
  Mathematica}, vol. 143, 01 2000.

\bibitem[Szarek(2000{\natexlab{c}})]{Szarek2000(3)}
------, ``Generic properties of learning systems,'' in \emph{Annales Polonici
  Mathematici}, vol.~73.\hskip 1em plus 0.5em minus 0.4em\relax Instytut
  Matematyczny Polskiej Akademii Nauk, 2000, pp. 93--103.

\bibitem[Horbacz and Szarek(2001)]{Szarek2001}
K.~Horbacz and T.~Szarek, ``Continuous iterated function systems on polish
  spaces,'' \emph{Bulletin of the Polish Academy of Sciences, Mathematics},
  vol.~49, 01 2001.

\bibitem[Szarek(2003{\natexlab{b}})]{Szarek2003(2)}
T.~Szarek, ``Invariant measures for non-expansive markov operators on polish
  spaces,'' \emph{Dissertationes Mathematicae}, vol. 415, pp. 1--62, 01 2003.

\bibitem[Ghosh et~al.(2019)Ghosh, Marecek, and Shorten]{Ramen2019}
R.~Ghosh, J.~Marecek, and R.~Shorten, ``Iterated piecewise-stationary random
  functions,'' \emph{arXiv preprint arXiv:1909.10093}, 2019.

\bibitem[Kifer(2012)]{Kifer2012}
Y.~Kifer, \emph{Ergodic theory of random transformations}.\hskip 1em plus 0.5em
  minus 0.4em\relax Springer Science \& Business Media, 2012, vol.~10.

\bibitem[Bhattacharya and Waymire(2009)]{Bhattacharya2009}
R.~N. Bhattacharya and E.~C. Waymire, \emph{Stochastic processes with
  applications}.\hskip 1em plus 0.5em minus 0.4em\relax SIAM, 2009.

\bibitem[Schur(1921)]{Schur1921}
\BIBentryALTinterwordspacing
J.~Schur, ``Über die gaußschen summen,'' \emph{Nachrichten von der
  Gesellschaft der Wissenschaften zu Göttingen, Mathematisch-Physikalische
  Klasse}, vol. 1921, pp. 147--153, 1921. [Online]. Available:
  \url{http://eudml.org/doc/59098}
\BIBentrySTDinterwordspacing

\bibitem[Graham and Lehmer(1976)]{Graham1976}
R.~L. Graham and D.~H. Lehmer, ``On the permanent of schur's matrix,''
  \emph{Journal of the Australian Mathematical Society}, vol.~21, no.~4, p.
  487–497, 1976.

\bibitem[Bof et~al.(2018)Bof, Carli, and Schenato]{Bof2018}
N.~Bof, R.~Carli, and L.~Schenato, ``Lyapunov theory for discrete time
  systems,'' \emph{arXiv preprint arXiv:1809.05289}, 2018.

\bibitem[Werner(2005)]{Werner2005}
I.~Werner, ``Contractive {M}arkov systems,'' \emph{Journal of the London
  Mathematical Society}, vol.~71, no.~1, pp. 236--258, 2005.

\bibitem[Werner(2004)]{Werner2004}
------, ``Ergodic theorem for contractive {M}arkov systems,''
  \emph{Nonlinearity}, vol.~17, no.~6, pp. 2303--2313, 2004.

\bibitem[Epperlein and Mareček(2017)]{8262886}
J.~Epperlein and J.~Mareček, ``Resource allocation with population dynamics,''
  in \emph{2017 55th Annual Allerton Conference on Communication, Control, and
  Computing (Allerton)}, 2017, pp. 1293--1300.

\bibitem[Kleywegt and Shao(2021)]{kleywegt2021optimizing}
A.~J. Kleywegt and H.~Shao, ``Optimizing pricing, repositioning, en-route time,
  and idle time in ride-hailing systems,'' \emph{arXiv preprint
  arXiv:2111.11551}, 2021.

\bibitem[Jiang et~al.(2021)Jiang, Kong, and Zhang]{Jiang2021}
\BIBentryALTinterwordspacing
Z.-Z. Jiang, G.~Kong, and Y.~Zhang, ``Making the most of your regret: Workers'
  relocation decisions in on-demand platforms,'' \emph{Manufacturing \& Service
  Operations Management}, vol.~23, no.~3, pp. 695--713, 2021. [Online].
  Available: \url{https://doi.org/10.1287/msom.2020.0916}
\BIBentrySTDinterwordspacing

\end{thebibliography}
\clearpage
\onecolumn

\newpage

\section*{Supplemental Material}
\vspace{0.75cm}
\setcounter{page}{1}

\section{An Overview of Notation}
\label{app:symbols}

The following notation is commonly used throughout the main manuscript and the supplemental material.

\begin{tabularx}{\linewidth}{ l | X }
\caption{Table of Notation}\\\toprule\endfirsthead
\toprule\endhead
\midrule\multicolumn{2}{l}{\itshape continues on next page}\\\midrule\endfoot
\bottomrule\endlastfoot
\textbf{Symbol} & \textbf{Meaning} \\\midrule
$\mathbb{N}$      & the set of natural numbers. \\
${\mathbb Z}$     & the set of all integers.\\
$\mathbb{Q}$      & the set of rational numbers. \\
$\mathbb{R}$      & the set of real numbers.  \\
$\mathbb{G}$        & a generic event. \\
$\mathbb{A}_i$      & the set of $i$th agent's actions. \\
$\mathbb{X}_i$      & a private state space of agent $i$, often $\R^{n_i}$. \\
$\mathbb{X}_S$      & Cartesian product of all agent's state space. \\
$\mathbb{X}^1_F$    & a space of internal states of the filter. \\
$\mathbb{X}^1_C$    & a space of internal states of the central controller. \\
$\mathbb{X}$        & combined state-space of the controller, the filter, and the agents. \\
$\mathbb{D}_i$      & set of possible resource demands of agent $i$. \\
${\mathbb O}_{\mathcal F}$& set of possible output values of the filter ${\mathcal F}$.\\
$\mathcal{B}(\mathbb X)$& Borel $\sigma$-algebra.  \\
$\mathcal M(\mathbb{X})$& a measure-space over $\mathbb{X}$. \\
$\mathcal M(\mathbb{X}^{\infty})$& a measure space over the path space.\\
${\mathcal E}$      & a real additive group. \\
${\mathcal C}_{\text{coup}}$& set of couplings.\\
$\mathcal{H}$       & a set used in the definition of asymptotic couplings. \\
${\mathcal{C}}^1$ & controller representing the central authority  for the first side of the two-sided market.\\
${\mathcal{C}}^2$ & controller representing the central authority  for the second side of the two-sided market. \\
${\mathcal F}^1$  & filter for the first side of the two-sided market.\\
${\mathcal F}^2$  & filter for the second side of the two-sided market.\\
${\mathcal{S}_1}^1, \ldots, {\mathcal{S}_N}^1$& systems modelling customers (e.g., seeking a ride). \\
${\mathcal{S}_1}^2, \ldots, {\mathcal{S}_N}^2$& systems modelling workers (e.g., drivers).\\
${\mathcal H}_{ij}$ & an output map. \\
$\mathcal{A}$       & augmented state transition matrix.\\
$\mathcal P_k$      & a measure-space-over-states-to-measure-space-over-states operator. \\
$\mathcal{P}_i^1$&Population 1 in Toy Example 1, in Figure \ref{f1}.\\
$\mathcal{P}_i^2$&Population 2 in Toy Example 1, in Figure \ref{f1}.\\
$\Tilde{\mathcal H}$& a block matrix used in \eqref{eq:int-descrptn}.\\
$\bbPi$             & the set of admissible broadcast control signals. \\
$\Gamma$            & a generic measure over the product of the two path spaces, potentially a coupling.\\
$A^1_c$             & a matrix used in the controller of Theorem \ref{thm:ITE-NLS-LCF}.\\
$A^2_c$             & a matrix used in the controller of Theorem \ref{thm:ITE-NLS-LCF}.\\
$B^1_c$             & a matrix used in the controller of Theorem \ref{thm:ITE-NLS-LCF}.\\
$A^2_f$             & a matrix used in the filter  of Theorem \ref{thm:ITE-NLS-LCF}.\\
$B^2_c$             & a matrix used in the controller  of Theorem \ref{thm:ITE-NLS-LCF}.\\
$B^1_f$             & a matrix used in the filter  of Theorem \ref{thm:ITE-NLS-LCF}. \\
$C^1_c $            & a matrix used in the controller  of Theorem \ref{thm:ITE-NLS-LCF}. \\
$B^2_f$             & a matrix used in the filter  of Theorem \ref{thm:ITE-NLS-LCF}. \\
$B_{i}$             & a matrix used in the agent dynamics of Theorem \ref{thm:ITE-NLS-LCF}. \\
$C^2_c $            & a matrix used in the controller  of Theorem \ref{thm:ITE-NLS-LCF}. \\
$C^1_f$             & a matrix used in the filter  of Theorem \ref{thm:ITE-NLS-LCF}.\\
$D^1_c$             & a matrix used in the controller in \eqref{eq:C1}.\\
$C^2_f$             & a matrix used in the filter in \eqref{eq:F2}.\\
$D^2_c$             & a matrix used in the controller in \eqref{eq:C2}.\\
$c^2_i$             & a vector used in the agent dynamics in \eqref{eq:S2}.\\
$c^1_i$             & a vector used in the agent dynamics in \eqref{eq:S1}. \\
$h_i$             & the number of output maps ${\mathcal H}_{ij}$. \\
$l_{ij}$          & Lipschitz constant for a transition map. \\
$l'_{ij}$         & Lipschitz constant for an output map. \\
$n$               & a dimension of a generic state space $\mathbb X$. \\
$n^1_c$           & dimension of the state of the controller for the first side of the two-sided market. \\
$n^1_f$           & dimension of the state of the filter  for the first side of the two-sided market. \\
$n^2_c$           & dimension of the state of the controller  for the second side of the two-sided market. \\
$n^2_f$           & dimension of the state of the filter  for the second side of the two-sided market. \\
$n^1_i$           & dimension of the state of $i$th agent's private state.  \\
$m^1_i$           & the number of possible actions of agent $i$. \\
$m$               & an upper bound on the number of possible actions of any agent. \\
$f_j$             & a generic map in generic iterated random functions. \\
$g$               & a generic function. \\
$r^1$             & the reference value; i.e., desired value of $y^1(k)$. \\
$\overline{r}^1_i$& $i$th agent's expected share of the resource over the long run. \\
$w_i$             & the number of state transition maps ${\mathcal W}_{ij}$. \\
$n^2_i$           & dimension of the state of $i$th agent's private state  for the second side of the two-sided market.  \\
$m^2_i$           & the number of possible actions of agent $i$. \\
$m$               & an upper bound on the number of possible actions of any agent. \\
$r^2$             & the reference value; i.e., desired value of $y^1(k)$. \\
$\overline{r}^1_i$& $i$th agent's expected share of the resource over the long run. \\
$w_i$             & the number of state transition maps ${\mathcal W}_{ij}$. \\
$w_{i,j}$         & a Borel map in a family of Borel map.\\
$p_{i,j}$    & a Borel measurable probability function in a family of Borel map.\\
$z$               & ${\mathcal Z}$-transform variable. \\
$x_c^1(k)$          & controller's internal state at time instant $k$. \\
$x_f^1$             & utilisation  of resource by $i^{th}$ agent at time instant $k$. \\
$y^1(k)$            & total resource utilisation at time instant $k$. \\
$\hat y^1(k)$       & value of $y^1(k)$ filtered by filter $\mathcal{F}^1$. \\
$p_j$               & a probability function of a generic iterated function system. \\
$p_{ij}$            & a probability function for the choice of agent $i$'s transition map. \\
$p'_{i\ell}$        & a probability function for the choice of agent $i$'s output map. \\
$\alpha$          & a constant used in the PI controller or its lag approximant.  \\
$\beta$           & a constant used in a lag controller.  \\
$\eta$            & a lower bound on the values of probability functions. \\
$\kappa$          & a constant used in the PI controller or its lag approximant.\\
$\xi$             &augmented state-vector. \\
$\Phi^{(1)}$        & projector from a measure over the product of the two path spaces to a single path space.\\
$\Phi^{(2)}$        & projector from a measure over the product of the two path spaces to a single path space.\\
$d^2_i$             & a vector used in the agent dynamics \eqref{eq:S2}. \\
$X(k)$              & element of a generic state space.  \\
$\{ X(k)\}_{k\in\N}$& a generic Markov chain.\\
$\beta_\ell$        & is built from all $b_{ij}$, ${d}_{ij}$, and other signals.\\
$N^1$               & the number of participants on one side of the market. \\
$N^2$               & the number of participants on the other side of the market. \\
$P(x,\mathbb{G})$   & a generic transition operator. \\
$P$ & a state-and-signal-to-state transition operator. \\

$e^1(k)$            & the error signal at time $k$; i.e., $\hat y^1(k) - r$. \\
$u^p$               & an external input, \eqref{eq:mult-side}.\\
$e^p$               & an external input, \eqref{eq:mult-side}.\\
$\Tilde{e}$         & a column vector of inputs, \eqref{eq:int-descrptn}. \\
$\Tilde{u}$         & a column vector of external inputs, \eqref{eq:int-descrptn}.\\
$H_{pq}$            & a matrix used in \eqref{eq:mult-side}.\\
$\pi^1(k)$          & the signal broadcast at time $k$. \\
$\lambda$           & initial state (distribution). \\
$\mu$               & a generic measure, usually on state space $\mathbb X$.  \\
$P_\lambda$         & a probability measure induced on the path space. \\
$\mathbf{1}$      & a compatible vector of ones. \\
\end{tabularx}


\section{Expanding the Interconnection Results Further}

Not all interconnection setups immediately fit into the frameworks presented in Sections \ref{sec:inter-two-ens} and \ref{sec:inter-lrgscl-ensmbl} of the main manuscript, respectively. In what follows, two toy examples are introduced, together with further interconnection results, so that the setups of the toy examples are accommodated.

\subsection{Toy Example 1}
\label{sec:toy1}

\begin{figure}[t]
\centering
\includegraphics[width=0.75\columnwidth]{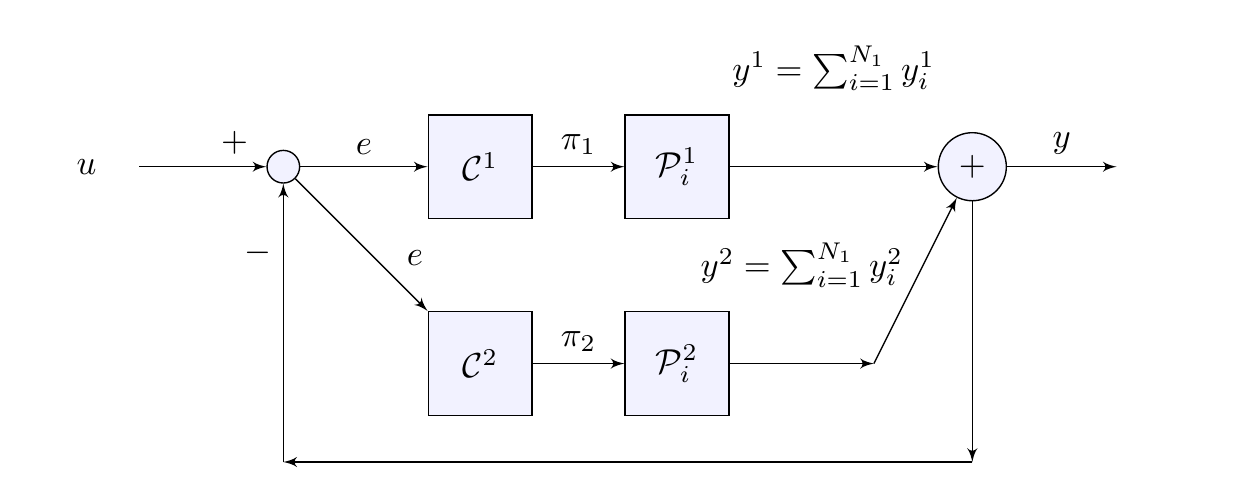}
\caption{The interconnection of two populations of drivers under consideration in Toy Example 1.}\label{f1}
\end{figure}

Consider the interconnection depicted in Figure \ref{f1}, where $\mathcal{P}_i^1$, $\mathcal{P}_i^2$, $\mathcal{C}^1$ and $\mathcal{C}^2$ are described by \eqref{eq:S1}, \eqref{eq:C1}, \eqref{eq:S2} and \eqref{eq:C2}, respectively, and (notation-wise) $e(k) := e^1(k) = e^2(k)$. In order to provide an interpretation for $\mathcal{P}_i^1$ and $\mathcal{P}_i^2$, let us suppose that each ensemble represents a different population of taxi drivers: $\mathcal{P}_i^1$ represents Population 1, and $\mathcal{P}_i^2$ represents Population 2. The drivers comprising Population 1 style themselves as being more open to offering a budget service and are thus relatively more likely to work (i.e., accept passengers) regardless of the fare price set by their hiring company. Meanwhile, the drivers that comprise Population 2 are more discerning in that they may not work (i.e., advertise rides) for lower-priced fares. The drivers' hiring companies set all fares.
Suppose that the demand signal, $u$, is the number of passengers in the community wanting rides. For simplicity, we will assume that a constant, ongoing demand throughout a time period of interest is maintained. The error signal is described by
$
    e(k) = u - (y^1(k) + y^2(k)),
$
where $y^1(k) = \sum^{N_1}_{i=1}y^1_i(k)$ and $y^2(k) = \sum^{N_2}_{i=1}y^2_i(k)$.

Let us assume that the controllers, $\mathcal{C}^1$ and $\mathcal{C}^2$, employed by each of the two hiring companies, are responsible for setting the prices of fares and operate similarly, the difference being that $\mathcal{C}^1$ updates every 40-time steps, while $\mathcal{C}^2$ is quicker, updating every 20-time steps. We will also assume that once a driver advertises a ride, the ride will be taken up by a waiting passenger, and the driver thus leaves the population for a new (or another returning) driver to take his or her place in the population. In other words, for simplicity, we will keep $N_1$ and $N_2$ fixed and constant.

The interconnection depicted in Figure \ref{f1} cannot be precisely categorised as one of the setups described in Theorems \ref{thm:ITE-NLS-LCF} or \ref{thm:lrg-scale-intr-cnnctn}. Therefore, we require the following proposition to proceed.

\begin{prop}\label{prop:toy_eg_1}
For the feedback system as in \ref{f1}, with $\mathcal{P}_i^1$, $\mathcal{P}_i^2$, $\mathcal{C}^1$ and $\mathcal{C}^2$ described by \eqref{eq:S1}, \eqref{eq:C1}, \eqref{eq:S2} and \eqref{eq:C2}, respectively. Furthermore, suppose that every agent (or driver) $i$ in a population has its state dynamics determined by the  stochastic difference equations as in \eqref{eq:S1} or \eqref{eq:S2}, where $A_i^1, A_i^2$ are Schur matrices, and $b_i^1, b_i^2$ and $d_i^1, d_i^2$ are selected, at every time instant according to probability functions which satisfies Dini continuity conditions in the manner of Theorem \ref{thm:ITE-NLS-LCF}, and that these probability functions are bounded below by scalars strictly more significant than 0. Then, for any pair of stable linear controllers $\mathcal{C}^1$ and $\mathcal{C}^2$ that are adaptable to the structure of the system,  the feedback loop converges in distribution to a unique invariant measure.

Proof: The proof follows in a manner similar to the proofs for Theorems \ref{thm:ITE-NLS-LCF} and \ref{thm:lrg-scale-intr-cnnctn}. That is, we define an augmented state
$
\xi(k) \defeq \left[
x^1(k)^T \,
x^2(k)^T \,
x_c^1(k)^T \,
x_c^2(k)^T
\right]^T,
$
where
$
x^1(k) \defeq \left[
x^1_1(k)^T \, x^1_2(k)^T \, \dots \ x^1_{N_1}(k)^T
\right]^T
$ and $
x^2(k) \defeq \left[
x^2_1(k)^T \, x^2_2(k)^T \, \dots \ x^2_{N_2}(k)^T
\right]^T,
$
whose dynamic behaviour is described by the difference equation
$
\xi(k+1) = {\mathcal W}_\ell(x) \defeq \mathcal{A}\xi(k) + \beta_\ell,   
$
where $\beta_\ell$ is built by linear combinations of $b^1_{ij}, b^2_{ij}$, the scalars ${d}^1_{ij}, {d}^2_{ij}$ and other signals, and
\begin{align}
\mathcal{A}:=
\begin{bmatrix}
\hat A^1 & 0 & 0 & 0\\
0 & \hat A^2 & 0 & 0\\
-B^1_c {\bf 1}^T \hat C^1 & -B^1_c {\bf 1}^T \hat C^2 & A_c^1 & 0 \\
-B^2_c {\bf 1}^T \hat C^1 & -B^2_c {\bf 1}^T \hat C^2 & 0 & A_c^2
\end{bmatrix}
\end{align}
$\mathbf{1}$ is a row vector consist of $1$ as each entry, ${\hat A}^1 \defeq \mathbf{diag}(A^1_i)$, ${\hat C}^1 \defeq
\mathbf{diag}(c^{1T}_i)$, ${\hat A}^2 \defeq \mathbf{diag}(A^2_i)$ and ${\hat C}^2 \defeq
\mathbf{diag}(c^{2T}_i)$. The proof then follows in a manner similar to the proof of Theorem \ref{thm:ITE-NLS-LCF}. $\blacksquare$
\end{prop}

\subsection{Toy Example 2}
\label{sec:toy2}
\begin{figure}[ht]
\centering
\includegraphics[width=0.8\columnwidth]{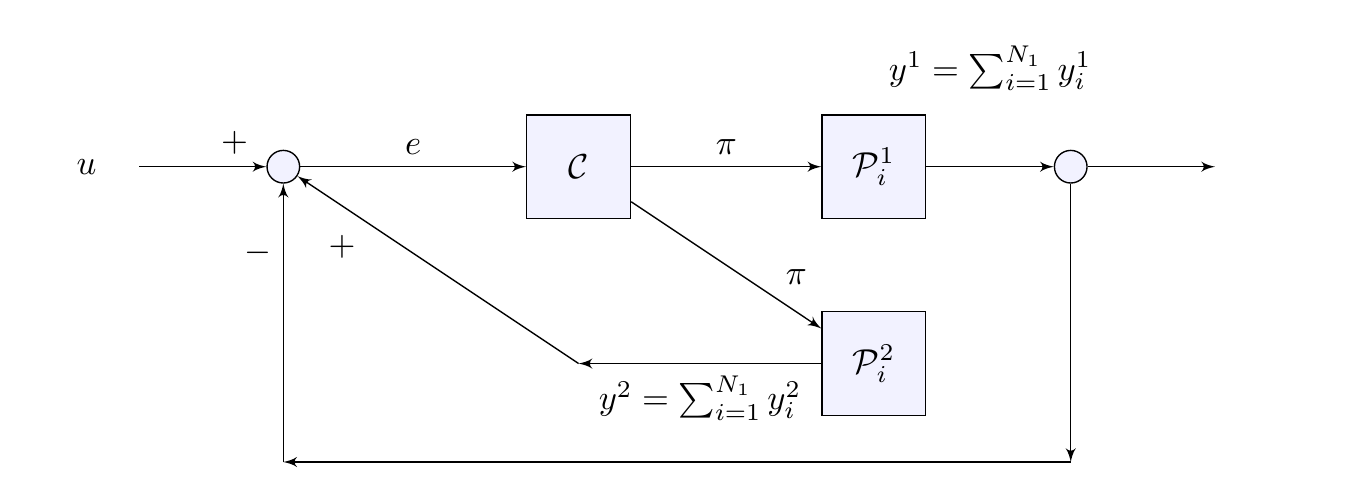}
 \caption{An interconnection with an elastic demand, as described in Toy Example 2.}\label{f8}
\end{figure}

For our second toy example, consider the interconnection depicted in Figure \ref{f8}, where $\mathcal{P}_i^1$ and $\mathcal{P}_i^2$ are described by \eqref{eq:S1} and \eqref{eq:S2}, respectively; and $\mathcal{C}$ is described by 
\begin{equation}\label{eq:C_toy2}
{\mathcal C} ~:~ \left\{ \begin{array}{rcl}
x_c(k+1) & = & A_c x_c(k) + B_c e(k), \vspace{0.1cm} \\
\pi(k) & = & C_c x_c(k) + D_c e(k).
\end{array} \right.
\end{equation}
The interpretation that we will give to $\mathcal{P}_i^1$ is that the ensemble again represents a population of taxi drivers. Similarly, the signal $u$ will again represent a fixed demand; that is, several passengers in the community seeking rides. The fixed demand remains constant throughout the experiment.

In Toy Example 2, however, we will provide a different interpretation for $\mathcal{P}_i^2$. Suppose that, in addition to the fixed demand, there is an elastic demand for rides that is influenced by the current price of a taxi fare as set by the company that hires the taxi drivers. This elastic demand is given by $y^2(k)$, where $y^2(k) = \sum^{N_2}_{i=1}y^2_i(k)$. As such, the total demand for rides is equal to $u + y^2(k)$. The error signal, $e$, is thus described by
$
    e(k) = u + y^2(k) - y^1(k),
$
where $y^1(k) = \sum^{N_1}_{i=1}y^1_i(k)$. As in Toy Example 1, $N_1$ and $N_2$ are assumed to be fixed and constant.

Similar to Toy Example 1, the interconnection depicted in Figure \ref{f8} cannot be precisely categorised as one of the setups described in Theorems \ref{thm:ITE-NLS-LCF} or \ref{thm:lrg-scale-intr-cnnctn}. Therefore, we again require a new proposition to proceed.

\begin{prop}\label{prop:toy_eg_2}
Let us note the Figure \ref{f8}, with $\mathcal{P}_i^1$, $\mathcal{P}_i^2$ and $\mathcal{C}$ described by \eqref{eq:S1}, \eqref{eq:S2} and \eqref{eq:C_toy2}, respectively. Furthermore, suppose that every agent $i$ in a population has its state dynamics determined by the stochastic difference equations given in \eqref{eq:S1} or \eqref{eq:S2}, where $A_i^1, A_i^2$ are Schur matrices, and $b_i^1, b_i^2$ and $d_i^1, d_i^2$ are chosen, at each time instant, with the probability functions which satisfies Dini conditions for continuity in the manner of Theorem \ref{thm:ITE-NLS-LCF}, and that these probability functions are bounded below by scalars strictly more significant than 0. Then, for any stable linear controller $\mathcal{C}$ adaptable to the configuration of the system, the feedback loop converges in distribution to a unique invariant measure.

Proof: The proof follows in a manner similar to the proofs for Theorems \ref{thm:ITE-NLS-LCF} and \ref{thm:lrg-scale-intr-cnnctn}, and Proposition \ref{prop:toy_eg_1}. We define an augmented state
$
\xi(k) \defeq \left[
x^1(k)^T \,
x^2(k)^T \,
x_c(k)^T
\right]^T,
$
where
$
x^1(k) \defeq \left[
x^1_1(k)^T \, x^1_2(k)^T \, \dots \ x^1_{N_1}(k)^T
\right]^T $ and $
x^2(k) \defeq \left[
x^2_1(k)^T \, x^2_2(k)^T \, \dots \ x^2_{N_2}(k)^T
\right]^T,
$
whose dynamic behaviour is described by the difference equation
$
\xi(k+1) = {\mathcal W}_\ell(x) \defeq \mathcal{A}\xi(k) + \beta_\ell,   
$
where $\beta_\ell$ is built from the linear combination of the
vectors $b^1_{ij}, b^2_{ij}$, the scalars ${d}^1_{ij}, {d}^2_{ij}$ and other signals, and
\begin{align}
\mathcal{A}:=
\begin{bmatrix}
\hat A^1 & 0 & 0\\
0 & \hat A^2 & 0\\
-B_c {\bf 1}^T \hat C^1 & B_c {\bf 1}^T \hat C^2 & A_c
\end{bmatrix}
\end{align}
where $\mathbf{1}$ is the vector of ones, ${\hat A}^1 \defeq \mathbf{diag}(A^1_i)$, ${\hat C}^1 \defeq
\mathbf{diag}(c^{1T}_i)$, ${\hat A}^2 \defeq \mathbf{diag}(A^2_i)$ and ${\hat C}^2 \defeq
\mathbf{diag}(c^{2T}_i)$. The proof then follows in a manner similar to the proof of Theorem \ref{thm:ITE-NLS-LCF}. $\blacksquare$
\end{prop}

\section{Numerical Illustrations and Discussion}
\label{sec:simulations}

To illustrate the meaning of our analytical results, we chose two intentionally simple examples,
which have been worked out in detail. We discuss how a variety of real-world settings could be
modelled in the proposed framework, too. 

\subsection{Toy Example 1}

\begin{figure}[ht]
\centering
\includegraphics[scale=1.5]{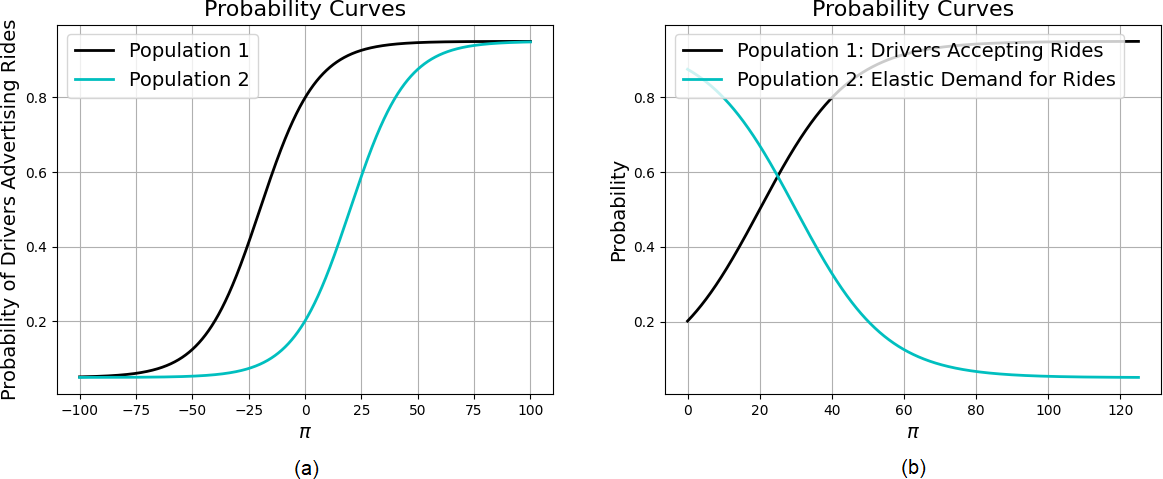}
\caption{Probability curves. (a) Probability of a driver advertising rides as a function of fare price $\pi$. (b) Probabilities as a function of fare price $\pi$.}\label{f2}
\end{figure}

First, let us consider the simple example where there are two populations of drivers as described in Section \ref{sec:toy1} above. 
The feedback loop is depicted in Figure \ref{f1}.
To corroborate Proposition \ref{prop:toy_eg_1}, we wish to demonstrate the convergence of the number of drivers available in distribution to a unique invariant measure.
We hence performed ten runs of a simulation with 1800 time steps in each run, with the following parameters: the sizes of Populations 1 and 2 were set at 50 and 100 drivers, respectively; the probability of drivers advertising rides at the beginning of each simulation was randomized using the Python function {\fontfamily{lmtt}\selectfont random.uniform(0,1)} (i.e., specifically, at the beginning of each simulation, the Python function was called once for the collective pool of drivers from Population 1, and then a second time for the collective pool of drivers belonging to Population 2); the probability of a driver advertising rides as a function of the currently set fare price, $\pi$, is indicated in Figure \ref{f2}(a); the fixed demand, $u$, was set to 120 for the duration of each simulation; and the controllers were described by $\pi(k)=\beta\pi(k-1)+\kappa[e(k)-\alpha e(k-1)]$, where $\alpha=-4.01$, $\beta=0.99$, and $\kappa=0.1$, noting that $\mathcal{C}^1$ updated every 40 time steps, while $\mathcal{C}^2$ updated every 20 time steps.

\begin{figure}[ht]
\centering
\includegraphics[scale=1.5]{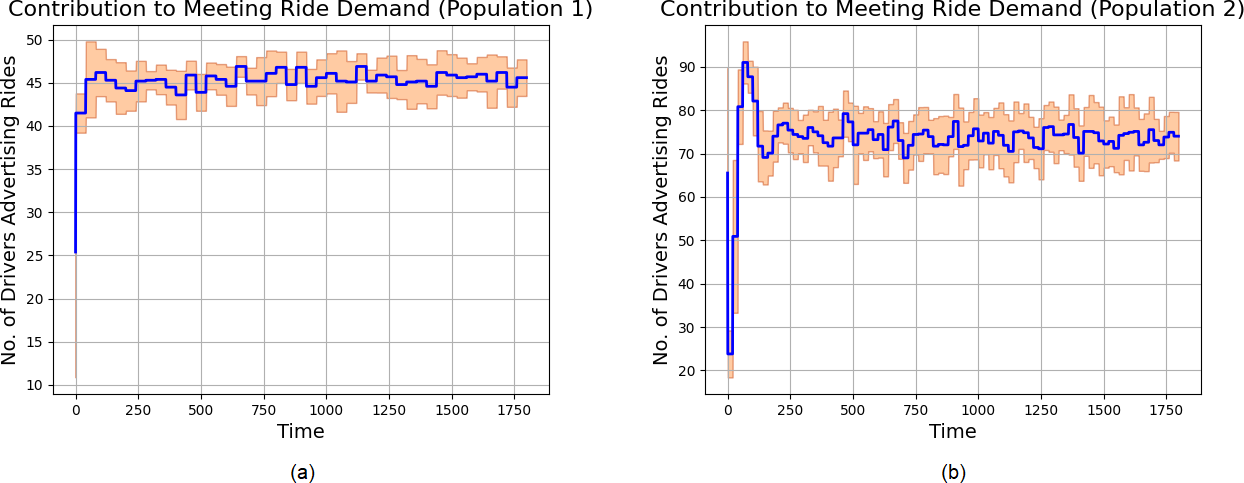}
\caption{Each population's contribution to meeting the demand for rides: (a) Population 1; (b) Population 2. The blue line indicates the mean number (from 10 simulation runs) of drivers advertising rides versus time, while the red shaded area indicates one standard deviation from the mean.}\label{f3}
\end{figure}

\begin{figure}[ht]
\centering
\includegraphics[scale=1.5]{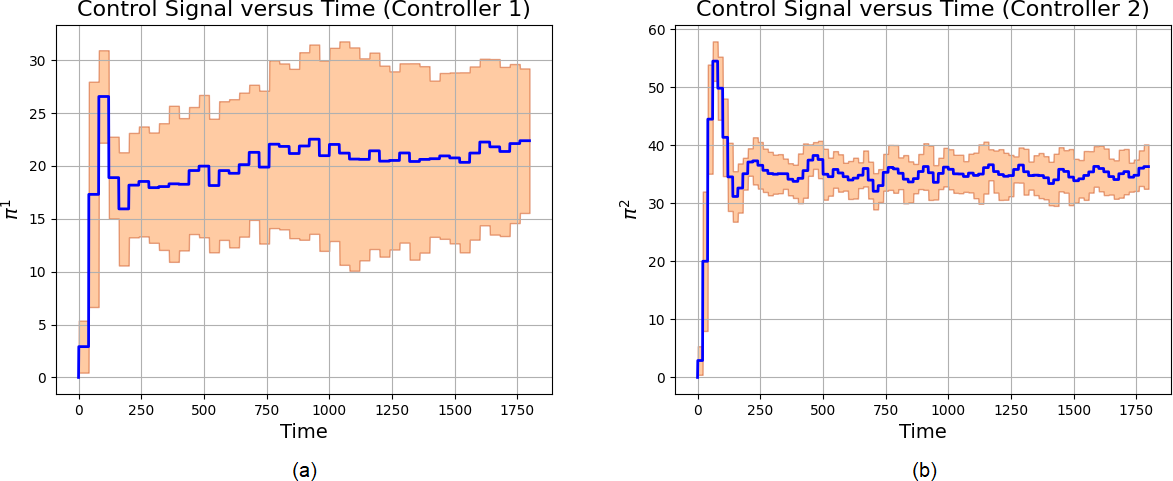}
\caption{(a) The evolution over time of the output from Controller 1, $\pi^1$, (i.e., the fare price set by the company to which Population 1 belongs). The blue line indicates the mean of $\pi^1$ versus time from 10 simulation runs, while the red shaded area indicates one standard deviation from the mean. The output of Controller 1 can be used in association with the black curve in Figure \ref{f2}(a). (b) The evolution over time of the output from Controller 2, $\pi^2$, (i.e., the fare price set by the company to which Population 2 belongs). The blue line indicates the mean of $\pi^2$ versus time from 10 simulation runs, while the red shaded area indicates one standard deviation from the mean. The output of Controller 2 can be used in association with the cyan curve in Figure \ref{f2}(a).}\label{f4}
\end{figure}

\begin{figure}[ht]
\centering
\includegraphics[scale=0.56]{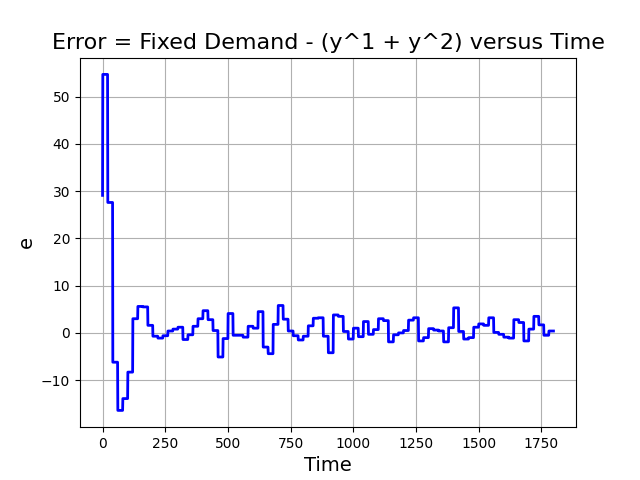}
\caption{The evolution of the error signal, $e$, over time. The blue line indicates the fixed demand, $u$, minus the sum of the means of $y^1$ and $y^2$ from 10 simulation runs.}\label{f7}
\end{figure}

The results from the experiment are as follows: Figures \ref{f3}(a) and \ref{f3}(b) show, on average, each driver population's contribution to meeting the demand for rides as time evolved; Figures \ref{f4}(a) and \ref{f4}(b) show, on average, the evolution of the control signals (i.e., the fare prices set by the two different taxi driver hiring companies) over time; and Figure \ref{f7} illustrates, on average, the evolution of the error signal, $e$, over time. From Figure \ref{f7}, it can be observed that the error signal evolution is contained near to zero; in other words, the demand for rides is sufficiently being met by the system.

\subsection{Toy Example 2}

Next, let us consider the other simple example, where the demand for rides is elastic; i.e., influenced by the current price of a fare, as introduced in Section~\ref{sec:toy2} above.
For the feedback model of the interconnection, see Figure \ref{f8}. The demand for rides as a function of the current price of a taxi fare is displayed in Figure \ref{f2}(b); see the cyan curve. As the price of fares decreases, the probability that passengers will want additional rides in the community (on top of the rides already wanted that comprise the fixed demand) increases. Meanwhile, as the price of fares decreases, the probability of a driver advertising rides also decreases; see the black curve.

To corroborate Proposition \ref{prop:toy_eg_2}, we wish to demonstrate the convergence of the number of drivers available in distribution to a unique invariant measure. We also wish to demonstrate the convergence of the elastic demand in distribution to a unique invariant measure. We perform ten runs of a simulation with 1800 time steps in each run, with the following parameters: the fixed demand, $u$, was set to 20 for the duration of each simulation; the sizes of Populations 1 and 2 were set to 60 drivers, and a value of 20 (i.e., the maximum additional demand on top of the fixed demand), respectively; the probability of drivers advertising rides at the beginning of each simulation was randomized using the Python function {\fontfamily{lmtt}\selectfont random.uniform(0,1)} (i.e., specifically, at the beginning of each simulation, the Python function was called once for the collective pool of drivers); similarly, at the beginning of each simulation, the Python function {\fontfamily{lmtt}\selectfont random.uniform(0,1)} was called once for the collective pool of agents representing elastic demand from Population 2, to determine the initial probability of an agent contributing to the elastic demand for that simulation run; and the controller was described by $\pi(k)=\beta\pi(k-1)+\kappa[e(k)-\alpha e(k-1)]$, where $\alpha=-4.01$, $\beta=0.99$, and $\kappa=0.1$. The controller produced an updated output every 20 time steps.

\begin{figure}[ht]
\centering
\includegraphics[scale=1.5]{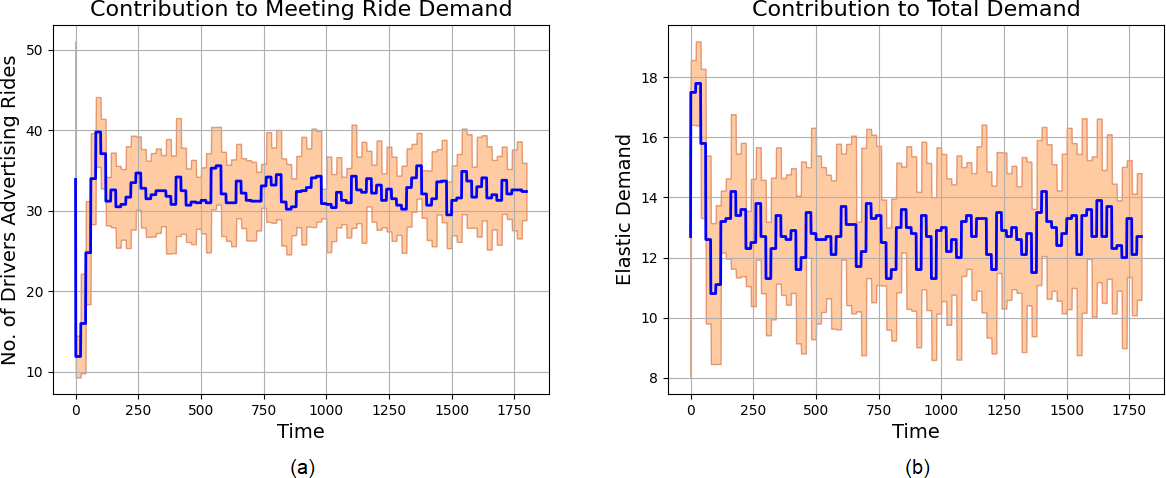}
\caption{(a) Population 1's contribution to meeting the demand for rides. The blue line indicates the mean number (from 10 simulation runs) of drivers advertising rides versus time, while the red shaded area indicates one standard deviation from the mean. (b) Population 2's contribution to the total demand for rides. The blue line indicates the mean number (from 10 simulation runs) of additional rides required (on top of the fixed demand) versus time, while the red shaded area indicates one standard deviation from the mean.}\label{f10}
\end{figure}

\begin{figure}[ht]
\centering
\includegraphics[scale=1.5]{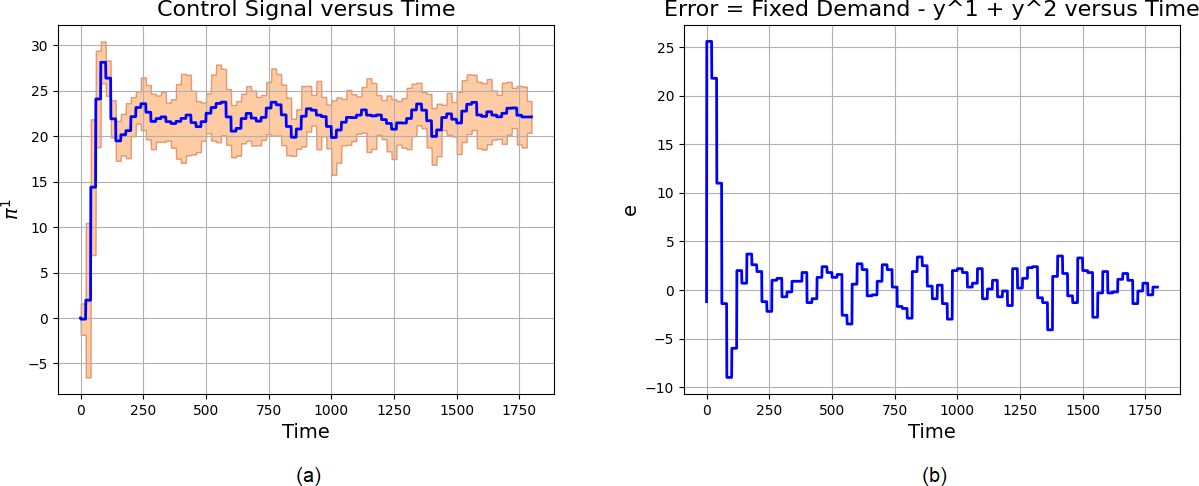}
\caption{(a) The evolution over time of the output from the Controller, $\pi^1$, (i.e., the fare price set by the company to which Population 1 belongs). The blue line indicates the mean of $\pi^1$ versus time from 10 simulation runs, while the red shaded area indicates one standard deviation from the mean. The output of the Controller can be used in association with the curves in Figure \ref{f2}(b). (b) The evolution of the error signal, $e$, over time. The blue line indicates the fixed demand, $u$, plus the mean of the elastic demand, $y^2$, minus the mean number of drivers advertising rides, $y^1$, from 10 simulation runs.}\label{f11}
\end{figure}

The results from the experiment are as follows: Figure \ref{f10}(a) shows, on average, the driver population's contribution to meeting the demand for rides as time evolved; Figure \ref{f10}(b) shows, on average, the elastic demand's contribution to the total demand for rides as time evolved; Figure \ref{f11}(a) shows, on average, the evolution of the control signal (i.e., the fare prices set by the taxi driver hiring company) over time; and Figure \ref{f11}(b) illustrates, on average, the evolution of the error signal, $e$, over time. From Figure \ref{f11}(b), it can be observed that the error signal evolution is contained near to zero; in other words, the demand for rides is sufficiently being met by the system.

\subsection{More Complicated Examples}

In the ride-hailing business of 
Uber Technologies, Inc., Lyft, Inc., or Didi Chuxing Technology Co.,
the feedback models are considerably more complicated than the two toy examples above; moreover, their precise nature is not known publicly. 
However, when one revisits the model of Figure~\ref{fig:system_ii}, one can imagine applying it to 
a model utilizing the following assumptions:
\begin{itemize}
    \item a two-sided market operated within one spatial region.
    For example, within New York City, Manhattan may be operated as one market. 
    Brooklyn would be another market, and the interconnection would not be modelled in our simulations, although Theorem \ref{thm:lrg-scale-intr-cnnctn} would apply to such interconnections.
    \item the number of installations of the app as a constant.
    We do allow for time-varying (or state-varying) probability of a particular customer seeking a ride at a particular time, cf. \cite{8262886}, which could model some systems coming online only at some later point.
    \item the number of registered driver-partners are constant during our operations.
    We do allow for time-varying (or state-varying) probability of a particular driver-partner accepting a matched ride at a particular time, cf. \cite{8262886}. A constant probability function can then model drivers going offline.
    \item a batched matching strategy that introduces a constant delay. 
    While a matching strategy could generally wait until there is a certain number of requests for rides, 
    we consider a matching strategy that waits for a constant interval, e.g., 5 seconds, between matching the available requests for rides to driver-partners. 
    \item only the driver surge pricing component of the price, which we let vary continuously. Furthermore, we assume that the drivers cannot set their surge prices, which is in line with the latest changes of the system (see Faiz Siddiqui: Where have all the Uber drivers gone? Washington Post, May 7, 2021, \url{https://www.washingtonpost.com/technology/2021/05/07/uber-lyft-drivers/}), 
    and we disregard relocation decisions \citep{kleywegt2021optimizing,Jiang2021}.
    In order to balance the spatial distribution of the driver-partners, a platform may control other components of the price too, and the driver surge pricing may take on integer values, or generally values from a discrete set, as in Theorem \ref{thm:discrete-range-space}, but we disregard these possibilities for the sake of the clarity of the presentation.
\end{itemize}

\begin{rem}
The additive surge pricing rounded to quarters of a dollar corresponds to, in control-engineering terms, a ``Deadband'' in the controller, in the sense that up to some magnitude of the change in the error signal, there is no change in the control signal.
This makes it impossible to apply Theorem~\ref{thm:intrcon-twoensmbl-nonlinrstate-lincontr&filtr}, but makes it possible to
apply Theorem \ref{thm:discrete-range-space}.
\end{rem}

In a straightforward scenario, one could assume that the drivers accept all matches made by the system based on some greedy first-come-first-served procedure. 
This would be reasonable because the drivers need to keep up their ``acceptance rate'' in order to receive matches
by the drivers (these policies are not made officially known but are widely observed, cf. \url{https://www.uberpeople.net/threads/ignore-vs-decline.310718/}). 
The corresponding ${\mathcal C^2}$ would be strikingly simple, would consider one request for a ride at each time, and would guarantee  $y^{2}(k) = \hat y^{1}(k)$ whenever there is sufficient capacity.
The filtered output $\hat y^2(k)$ could consist of the proportion of empty cars on the road, with the reference signal $u^1(k)$ ranging between 10-15\%. This suggests that there should be some empty cars on the road, but not too many,
as implemented by providing the error signal \eqref{eq:e1} to the controller ${\mathcal C^1}$.
The controller ${\mathcal C^1}$ suggests prices $\pi^1$ while considering additive driver surge pricing based on 
an inner state of the controller $x_c^1(k)$. 
From the results of \cite[Section 3.1]{Fioravanti2019}, it is easy to see that a PI controller may not be suitable for  ${\mathcal C^1}$. That is, for some constants $K_p, K_i$:
\begin{align}
\pi^1_{\text{PI}}(k+1)&=\left( K_p e^1(k)+K_i (x_c^1 (k)+e^1(k) \right) \label{pi1PI} \\ 
\pi^1_{\text{Lag}}(k+1)&=\left( K_p e^1(k)+K_i (0.99 x_c^1 (k)+e^1(k) \right) \label{pi1lag}
\end{align}
The PI controller \eqref{pi1PI} destroys the ergodicity, while its lag approximant \eqref{pi1lag} allows for the unique ergodicity. 
The customers from the set $\{{\mathcal{S}_i}^1\}_{i=1}^{N}$ respond to the signal $\pi^1$ by issuing or not issuing requests from the set $\{y^{1}_i(k)\}_{i=1}^{N}$ for rides at time $k$, based on some internal state of the customers $\{x^{1}_i(k)\}_{i=1}^{N}$ at time $k$, which are not directly observable.

\end{document}